\documentclass[11pt,a4paper]{article}
\usepackage{ifthen,latexsym,amssymb,amsmath,bbm,fixmath,stmaryrd}
\usepackage[nobysame,initials]{amsrefs}
\usepackage[shortlabels]{enumitem}
\usepackage{listings} 
\usepackage{hyperref}
\hypersetup{
    colorlinks=true,
    linkcolor=blue,
    filecolor=magenta,      
    urlcolor=cyan,
}

\usepackage[noabbrev]{cleveref}

\usepackage{algorithm2e}
\usepackage{tabularx}
\usepackage{float}
\usepackage{etoolbox}

\newcommand{\beq}[1]{\begin{equation}\label{#1}}
  \newcommand{\eeq}{\end{equation}}

\newcommand{\e}{\varepsilon}

\newcommand{\C}[1]{{\protect\mathcal{#1}}}
\newcommand{\B}[1]{{\bf #1}}
\newcommand{\I}[1]{{\mathbbm #1}}
\renewcommand{\O}[1]{\overline{#1}}

\newcommand{\V}[1]{\mathbold{#1}}

\newtheorem{definition}{Definition} [section]
\newtheorem{theorem}[definition]{Theorem}
\newtheorem{lemma}[definition]{Lemma}

\newtheorem{conjecture}[definition]{Conjecture}
\newtheorem{claim}[definition]{Claim}

\newtheorem{observation}[definition]{Observation}

\usepackage{changepage}

\usepackage{pgf,tikz}
\usepackage{pdflscape}

\newcommand{\hide}[1]{}

\newcommand{\bpf}[1][Proof.]{\smallskip\noindent{\it #1} }
\newcommand{\qed}{\nolinebreak\mbox{\hspace{5 true pt}%
  \rule[-0.85 true pt]{3.9 true pt}{8.1 true pt}}}
\newcommand{\epf}{\qed \medskip}

\newcommand{\formatnum}[2]{
  \pgfmathprintnumber[fixed, precision=\ifstrempty{#1}{10}{#1}]{#2}...
}

\begin{document}


\newcommand{\PST}[1]{\ifthenelse{\equal{#1}{}}
{\cite{PikhurkoSliacanTyros19}}
{\cite{PikhurkoSliacanTyros19}*{#1}}}

\newcommand{\LPSS}[1]{\ifthenelse{\equal{#1}{}}
{\cite{LiuPikhurkoSharifzadehStaden23}}
{\cite{LiuPikhurkoSharifzadehStaden23}*{#1}}}

\newcommand{\EL}[1]{\ifthenelse{\equal{#1}{}}
{\cite{EvenzoharLinial15}}
{\cite{EvenzoharLinial15}*{#1}}}

\newcommand{\BGHKS}[1]{\ifthenelse{\equal{#1}{}}
{\cite{BasitGranetHorsleyKundgenStaden25x}}
{\cite{BasitGranetHorsleyKundgenStaden25x}*{#1}}}

\newcommand{\obj}[2][H]{\lambda_{#1}\ifstrempty{#2}{}{(#2)}}
\newcommand{\Obj}[2][H]{\Lambda_{#1}\ifstrempty{#2}{}{(#2)}}

\newcommand{\N}{N}
\newcommand{\p}{p}
\renewcommand{\P}{P}
\newcommand{\aut}{\mathrm{aut}}
\newcommand{\tinj}{t}
\renewcommand{\sqcup}{+}
\renewcommand{\mid}{:}

\newcommand{\blow}[2]{\ifstrempty{#2}
{#1{()}}
{#1(#2)}
}

\newcommand{\dedit}{\delta_{\mathrm{edit}}}

\renewcommand{\ge}{\geqslant}
\renewcommand{\geq}{\geqslant}
\renewcommand{\le}{\leqslant}
\renewcommand{\leq}{\leqslant}
\renewcommand{\preceq}{\preccurlyeq}
\renewcommand{\succeq}{\succcurlyeq}
\newcommand{\SectionH}[1]{\subsection{Graph $H_{#1}$}\label{se:#1}}
\newcommand{\Turan}[2]{T^{#1}_{#2}}


\newcommand{\construction}[2]{#1}
\newcommand{\PerfectStability}[2]{#1 &#2}
\newcommand{\MultiPartite}[1]{T_{#1}}

\tikzset{
  vertex/.style={
    circle,
    draw,
    fill=black,
    minimum size=3pt,
    inner sep=0pt
  },
  edge/.style={},
  edgeR/.style={red,line width=1.2pt},
  edgeB/.style={blue,line width=1.2pt}, 
  vertexlabel/.style={
    font=\scriptsize,
    inner sep=1pt
  }
}

\newcommand{\DrawGraph}[2]{
  \begin{tikzpicture}
    \def\graphradius{10pt} 
    \foreach \i in {0, ..., \the\numexpr #1-1 \relax} {
      \pgfmathsetmacro{\angle}{90 - \i * (360/#1)}
      \node[vertex] (v\i) at (\angle:\graphradius) {};
    }
    \foreach \u / \v in {#2} {
      \draw[edge] (v\u) -- (v\v);
    }
  \end{tikzpicture}
}

\usetikzlibrary{calc}

\newcommand{\DrawSemiGraph}[3]{
  \begin{tikzpicture}[baseline={($(current bounding box.center)-(0pt,3pt)$)}]
    \def\graphradius{7pt} 
    \foreach \i in {0, ..., \the\numexpr #1-1 \relax} {
      \pgfmathsetmacro{\angle}{90 - \i * (360/#1)}
      \node[vertex] (v\i) at (\angle:\graphradius) {};
    }
    \foreach \u / \v in {#2} {
      \draw[edgeR] (v\u) -- (v\v);
    }
    \foreach \u / \v in {#3} {
      \draw[edgeB] (v\u) -- (v\v);
    }
  \end{tikzpicture} %
}

\newcommand{\DrawSemiGraphLabeled}[4]{%
  \begin{tikzpicture}[baseline={($(current bounding box.center)-(0pt,3pt)$)}]
    \def\graphradius{9pt}
    \def\labelradius{15pt}
    \foreach \i in {0, ..., \the\numexpr #1-1 \relax} {
      \pgfmathsetmacro{\angle}{90 - \i * (360/#1)}
      \node[vertex] (v\i) at (\angle:\graphradius) {};
    }
    \foreach \lab [count=\i from 0] in {#4} {
      \pgfmathsetmacro{\angle}{90 - \i * (360/#1)}
      \node[vertexlabel] at (\angle:\labelradius) {\lab};
    }
    \foreach \u / \v in {#2} {
      \draw[edgeR] (v\u) -- (v\v);
    }
    \foreach \u / \v in {#3} {
      \draw[edgeB] (v\u) -- (v\v);
    }
  \end{tikzpicture}%
}

\newcommand{\loopB}{L}

\newcommand{\TheoremTemplate}[4]{It holds that $\obj[#1]{}=#4$, the problem is perfectly $B$-stable with $B = #2$, and the unique maximizer for $\obj[#1]{}\left(\blow{B}{\V a}\right)$ is $\V a = \left( #3 \right)$.
}

\newcommand{\fallingfact}[1]{ \, \underline{#1} \, }

\title{Semi-inducibility of 4-vertex graphs}

\author{Levente Bodn\'ar\thanks{Email: bodnalev@gmail.com
}\ \mbox{ and }Oleg Pikhurko\thanks{Email: O.Pikhurko@warwick.ac.uk}\\
Mathematics Institute and DIMAP\\
University of Warwick\\
Coventry CV4 7AL, UK}

\maketitle

\begin{abstract}
For a graph $H$ whose edges are coloured blue or red, the \emph{$H$-semi-inducibility problem} asks for the maximum, over all graphs $G$ of given order $n$, of the number of injections from the vertex set of $H$ into the vertex set of $G$ that send red (resp.\ blue) edges of $H$ to edges (resp.\ non-edges) of~$G$. We consider all possible 4-vertex non-complete graphs $H$ and essentially resolve all remaining cases except when $H$ is the 3-edge path 
coloured blue-blue-red in this order (or is equivalent to this case).
Some of our proofs are computer-generated, using the flag algebra method of Razborov.
\end{abstract}

\section{Introduction}

Basit, Granet, Horsley, K\"undgen and Staden~\cite{BasitGranetHorsleyKundgenStaden25x} introduced the following problem for a given blue-red edge-coloured graph $H$ with $\kappa$ vertices. For an (uncoloured) graph $G$, let $\Obj[H]{G}$ be the number of \emph{embeddings} of $H$ into $G$, that is, injections $V(H)\to V(G)$ that map red (resp.\ blue) edges of $H$ to edges (resp.\ non-edges) of $G$. If $G$ has $n\ge \kappa$ vertices then we also consider the normalised version 
\[
 \obj[H]{G}:=\frac{\Obj[H]{G}}{n^{\fallingfact{\kappa}}},
 \]
  where $n^{\fallingfact{\kappa}}:=n(n-1)\dots(n-\kappa+1)$ is the $\kappa$-th falling power of $n$. Thus $\obj[H]{G}$ is the probability that a random injection $V(H)\to V(G)$ is an embedding of $H$ into~$G$.
The \emph{$H$-semi-inducibility problem} asks for $\Obj[H]{n}$, the maximum of $\Obj[H]{G}$ over all graphs $G$ of order $n$ (i.e.\ with $n$ vertices) or, equivalently, for its normalised version $\obj[H]{n}:=\Obj[H]{n}/n^{\fallingfact{\kappa}}$ for $n\ge \kappa$.  It is easy to show (see \cite[Proposition~3.1]{BasitGranetHorsleyKundgenStaden25x}) that the limit \[
\obj[H]{}:=\lim_{n\to\infty}\obj[H]{n}
\]
 exists. We call $\obj[H]{}$ the \emph{semi-inducibility constant} of $H$.

Note that our parameter $\Obj[H]{G}$ is slightly different from the one used in~\cite{BasitGranetHorsleyKundgenStaden25x} where the authors count the pairs of sets of edges and non-edges in $G$ which, when coloured red and blue respectively, form a copy of~$H$. 
Of course, each such pair corresponds to the same number of embeddings of $H$ into $G$ (which is the number of colour-preserving automorphisms of $H$), so it is easy to change between these two parameters. We use $\Obj[H]{G}$ as it admits the natural probabilistic version $\obj[H]{G}$, which is convenient in some arguments.

If $H$ is a colouring of the complete graph on $V(H)$ then the $H$-semi-inducibility problem amounts to
the well-known and actively studied inducibility problem (that was introduced by Pippenger and Golumbic~\cite{PippengerGolumbic75} in 1975) for the red subgraph of $H$. This was the main motivation in~\cite{BasitGranetHorsleyKundgenStaden25x} for introducing this question.

Various results on the semi-inducibility problem have been proved in~\cite{BasitGranetHorsleyKundgenStaden25x,BodnarPikhurko25ind,ChenClemenNoel25x,ChenNoel25x,BaloghLidicky26semiind}; we refer the reader to these papers for further details.

Here we systematically analyse the case when $H$ has 4 vertices (as all cases with $v(H)\le 3$ are trivial). We exclude the cases when $H$ is complete (as then we obtain the well-studied graph inducibility problem), or $H$ has no edges of some colour (as then trivially $\obj[H]{n}=1$), or $H$ has an isolated vertex (as this reduces to a 3-vertex case). Since the extremal function $\obj[H]{n}$ does not change when we exchange the two colours in $H$, it is enough to consider only one graph from each such pair. With the above conventions and exclusions, there are, up to isomorphism, 18 possible 4-vertex 2-coloured graphs $H$ to consider (some of which have been previously resolved) and all of these are listed in Table~\ref{ta:semiind}. 

In brief, our contribution is  to resolve all but one of the remaining cases in the table. 

This project initially started as our attempt to see what the flag algebra method of Razborov would give for this new extremal problem. The method turned out to be highly successful in proving upper bounds on the semi-inducibility constant, giving the correct value of it for all but one non-complete 4-vertex graph $H$.
These upper bounds on $\obj{}{}$  are recorded in Table~\ref{ta:semiind}. The table also gives a construction giving the matching lower bound on $\obj[H_i]{}$, where we denote the complete partite graph with $m$ parts of sizes $n_0,\dots,n_{m-1}$ by $\MultiPartite{n_0,\,\dots,\,n_{m-1}}$ (and we omit ceiling and floor signs when these are not essential). Although flag algebras work with a limit version of the problem as $n\to\infty$ and directly give asymptotic results only, we complemented these calculations with various stability-type arguments to get more precise structural results about (almost) extremal graphs.

In six of the resolved cases $H_i$, namely for $i=4,5,6,9,10,11$, we proved that, for all sufficiently large $n$, every extremal graph $G$ is a blow-up of some  small $m$-clique $K_m$ (with $m=m(i)$ being at most $5$ in each of these cases). Thus one can consider the problem as combinatorially solved for large $n$: the determination of the exact value of $\Obj[H_i]{n}$ (resp.\ the set of extremal graphs) amounts to maximising an explicit polynomial over non-negative integers $n_0,\dots,n_{m-1}$ that sum up to~$n$ (resp.\ describing all optimal $m$-tuples $(n_0,\dots,n_{m-1})$). In fact, we prove in each of these cases that there is a unique $m$-tuple of non-negative reals $x_0\ge\dots\ge x_{m-1}$ summing up to $1$ such that every almost extremal graph of order $n\to\infty$ can be made into  $\MultiPartite{x_0n,\,\dots,\,x_{m-1}n}$ by changing at most $o(n^2)$ adjacencies. Thus all almost extremal graphs are close to each other in the edit distance.

The remaining studied cases, namely $H_0$, $H_1$ and $H_{15}$, although exhibiting some variety of extremal structures, are rather straightforward and we have been able to obtain an essentially full picture in these cases.

Finally, let us discuss the two cases where we could not determine the semi-inducibility constant via the plain flag algebra method.

One of these two cases is $H_5$ which is the red 3-edge star $\MultiPartite{3,\,1}$ with one blue edge added. Here a lower bound
\beq{eq:Index5}
\obj[H_5]{}\ge \frac{-171 \sqrt{57} + 7879}{43904}
\eeq
 follows by taking the complete 5-partite graphs with part sizes $(\alpha n, \alpha n, \alpha n, \alpha n, (1-4\alpha) n)$, where $\alpha:=(13-\sqrt{57})/{56}$.
 The plain flag algebra method seems to return the same numerical value as in the right-hand side of~\eqref{eq:Index5} but we have not been able to convert it to an exact-arithmetic proof.
 Instead, we observe  (in Lemma~\ref{lm:CompletePartite} here) that the symmetrisation method of Brown and Sidorenko~\cite{BrownSidorenko94} can be applied to this case to show that, for every $n\ge 4$, at least one extremal graph for $\obj[H_5]{n}$ is complete partite.%
Working with complete partite graphs greatly reduces the scale of flag algebra calculations and we were able to produce an exact-arithmetic certificate with entries in the field $\I Q(\alpha)$, proving that~\eqref{eq:Index5} is equality.

The only case in Table~\ref{ta:semiind} with $\obj[H]{}$ still unknown is $H_3$, which is the 3-edge path coloured blue-blue-red, say visiting vertices $0,1,2,3$ in this order. Here we have the following construction $G$ of order $n\to\infty$: take a disjoint union of an almost $\gamma \beta n$-regular graph on a vertex set $B$ with $|B|=(\beta+o(1)) n$ and a clique on a set $A$ with $|A|=n-|B|$.
The number of embeddings of $H_3$ into $G$ can be approximated as follows: we sum over all maps $f$ of  $1,2\in V(H)$ to non-adjacent vertices $f(1),f(2)\in V(G)$ of the product of the non-degree of $f(1)$ and the degree of $f(2)$, as these are the number of choices for $f(0)$ and $f(3)$ within factor $1+o(1)$.
 Thus $\obj{G}$, within additive $o(1)$, is 
\begin{equation}\label{eq:H3p}
 p(\beta,\gamma):=
 \beta(1-\beta) \cdot (1-\gamma \beta )(1-\beta) + 
 (1-\beta)\beta \cdot \gamma \beta^2 + 
 \beta^2 (1 - \gamma) \cdot (1-\beta \gamma) \beta\gamma, 
\end{equation}
 where the 3 products give the density of embeddings $f:H\to G$ with the pair $(f(1),f(2))$ belonging to $B\times A$, $A\times B$ and $B\times B$ respectively.
Routine optimisation  shows that the maximum of $p$ over $0\le \gamma, \beta\le 1$ is attained when $\beta = \formatnum{}{0.39829918}$ and $\gamma = \formatnum{}{0.28158008}$
satisfy
\begin{eqnarray}
    \beta^5 - \frac{10}{3}\,\beta^4   + \frac{2521}{576}\,\beta^3 - \frac{407}{144}\,\beta^2 + \frac{43}{48}\,\beta - \frac19&=&0,\label{eq:H3beta} \\
    \gamma^4   - \frac{37}{16}\,\gamma^3 + \frac{57}{32}\,\gamma^2 - \frac{9}{16}\,\gamma + \frac{1}{16}&=&0,\label{eq:H3gamma}
\end{eqnarray}
 giving  $\formatnum{}{0.150083407311578}$ as the maximal value. Here, the floating-point value returned by flag algebras seems to match this lower bound. Unfortunately, we have not been able to round it, so we leave this as an open problem.

\begin{conjecture}\label{cj:H3} For the 4-vertex graph $H_3$, which has blue edges $\{0,1\}$ and $\{1,2\}$, and a red edge $\{2,3\}$, it holds that $\obj[H_3]{}=p(\beta,\gamma)$, where $p$ is defined in~\eqref{eq:H3p} while $\beta$ and $\gamma$ are the appropriate roots of~\eqref{eq:H3beta} and~\eqref{eq:H3gamma} respectively.\end{conjecture}

\newcommand{\TableHelper}[8]{
#1 & \DrawSemiGraph{#2}{#3}{#4} & $#5$ & $#6$ & #7 \\ \hline
}

\renewcommand{\arraystretch}{1.5}

\begin{table}[h!]
\begin{center}
\begin{tabular}{c|c|c|c|c}
\hline
$i$ & $H_i$ & $\obj[H_i]{}$ & Construction & Reference  \\
\hline
\hline 

\TableHelper{0}{4}{0/1}{2/3}{\frac{1}{4}}{\mbox{Graphs of edge density $\frac12$}}{Section~\ref{se:H0}}{Any edge density $1/2$ works}
\TableHelper{1}{4}{0/1}{0/2,0/3}{\frac{4}{27}}{\mbox{$\frac n3$-regular graphs}}{Section~\ref{se:H1}}{Any 1/3 degree regular graph works}
\TableHelper{2}{4}{0/1}{1/2,0/3}{\frac{4}{27}}{\mbox{$\frac n3$-regular graphs}}{\cite[Theorem~6.2]{BodnarPikhurko25ind},\\ &&&&\cite[Theorem~1.2]{ChenClemenNoel25x}
}{Quasirandom 1/3 we proved}
\TableHelper{3}{4}{1/2}{0/1,0/3}{\le \formatnum{}{0.1500834091519}}{\ge \formatnum{}{0.150083407311578}}{Conjecture~\ref{cj:H3}}{-}
\TableHelper{4}{4}{0/2,0/3,1/2}{0/1}{\frac{4}{27}}{T_{{n}/{3},\, {n}/{3},\, {n}/{3}}}{Section~\ref{se:H4}}{Condition 3 fails}
\TableHelper{5}{4}{0/1,0/2,0/3}{2/3}{\frac{-171 \sqrt{57} + 7879}{43904}}{T_{\alpha n,\, \alpha n,\, \alpha n,\, \alpha n,\, (1-4\alpha) n}}{Section~\ref{se:H5}}{Multipartite stable}
\TableHelper{6}{4}{1/2}{0/2,0/3,2/3}{\frac{1}{8}}{T_{{n}/{2},\, {n}/{2}}}{Section~\ref{se:H6}}{Condition 3 fails}
\TableHelper{7}{4}{0/3}{3/2,1/2,1/0}{\frac{27}{256}}{\mbox{$\frac{1}{4}$-quasirandom graphs}}{\BGHKS{Theorem~1.9}}{Quasirandom?}
\TableHelper{8}{4}{3/2,0/3,1/2,1/0}{0/2}{\frac{4}{27}}{T_{{n}/{3},\, {n}/{3},\, {n}/{3}}}{\BGHKS{Theorem 8.1}}{Condition 3 fails}
\TableHelper{9}{4}{0/1,0/2,1/2,2/3}{0/3}{\frac{12}{125}}{T_{{n}/{5},\, {n}/{5},\, {n}/{5},\, {n}/{5},\, {n}/{5}}}{Section~\ref{se:computer}}{Stable}
\TableHelper{10}{4}{0/1,1/2}{0/2,0/3}{\frac{1}{8}}{T_{{n}/{2},\, {n}/{2}}}{Section~\ref{se:computer}}{Stable}
\TableHelper{11}{4}{0/2,0/3}{0/1,2/3}{\frac{1}{8}}{T_{{n}/{2},\, {n}/{2}}}{Section~\ref{se:computer}}{Stable}
\TableHelper{12}{4}{0/1,0/3}{1/2,2/3}{\frac{1}{8}}{T_{{n}/{2},\, {n}/{2}}}{\BGHKS{Theorem~1.7}}{Condition 3 fails}
\TableHelper{13}{4}{0/1,2/3}{0/3,1/2}{\frac{1}{8}}{T_{{n}/{2},\, {n}/{2}}}{\BGHKS{Theorem~1.6}}{Condition 3 fails}
\TableHelper{14}{4}{0/2,1/2,2/3}{0/1,0/3}{\frac{1}{8}}{T_{{n}/{2},\, {n}/{2}}}{\BGHKS{Section~8.1}}{Stable}
\TableHelper{15}{4}{0/2,0/3}{0/1,1/2,2/3}{\frac{1}{27}}{\mbox{Triangle-free $\frac{n}{3}$-regular graphs}}{Section~\ref{se:H15}}{Quasirandom?}
\TableHelper{16}{4}{0/3,2/3}{0/2,0/1,1/2}{\frac{1}{8}}{T_{{n}/{2},\, {n}/{2}}}{\BGHKS{Section~8.1}}{Stable}
\TableHelper{17}{4}{0/2,0/3,1/2}{0/1,2/3}{\frac{1}{8}}{T_{{n}/{2},\, {n}/{2}}}{\BGHKS{Section~8.1}}{Stable}

\end{tabular}
\caption{The values of the $H$-semi-inducibility constant for $4$-vertex non-complete graphs $H$.}
\label{ta:semiind}
\end{center}
\end{table}

\medskip
This paper is organised as follows. Section~\ref{se:Prelim} contains various definitions and auxiliary results needed in the paper. Then we discuss each individual graph $H_i$ in Section~\ref{se:proofs} (with Table~\ref{ta:semiind} also containing references to the corresponding sub-sections of Section~\ref{se:proofs}). First, we deal with those cases where we were able to prove the results without using flag algebras, which are the graphs $H_0$, $H_1$, $H_4$ and $H_{15}$. Among the remaining cases, we leave $H_5$ until the very end since it requires recalling some basic results from the theory of limits of complete partite graphs.

We provide a computer-generated certificate (and the scripts that can be used to generate or verify it) for each upper bound stated in Table~\ref{ta:semiind}, including also the cases that have been resolved elsewhere. On the other hand, when dealing with individual cases in Section~\ref{se:proofs}, we 
replaced computer-generated claims by human-readable proofs, whenever we could find a proof that had reasonable length and gave a better understanding of the problem.

\section{Preliminaries}\label{se:Prelim}

In this section, we present some definitions and auxiliary results that will be needed in the paper.


Let $\I R$ denote the set of reals. Let $\I N$ denote the set of non-negative integers and, for $n\in\I N$, we define $[n]:=\{0,\dots,n-1\}$. Note that we start indexing from $0$, merely to be consistent with the same convention as in our code. 
If the meaning is clear, we may abbreviate an unordered pair $\{u,w\}$ as $uw$, including the case when $u$ and $w$ are single-digit numbers.
For a set $X$ and an integer $\kappa\ge 0$, the set of all $\kappa$-subsets of $X$ is denoted by $\binom{X}{\kappa}$. Recall that $n^{\fallingfact{\kappa}}:=\prod_{i=0}^{\kappa-1}(n-i)$ denotes the $\kappa$-th falling power of~$n$. We may omit ceiling/floor signs when they are not essential.

For $n,s\in \I N$, let $\pi_s(n)$ be the maximum product of $s$ non-negative integers that sum to $n$. 
  Clearly, this is attained if and only if the $s$ summands are nearly equal. So we can equivalently define
  \begin{equation}\label{eq:pi}
  \pi_s(n):=\prod_{i=0}^{s-1} \left\lfloor \frac{n+i}s\right\rfloor.
  \end{equation}
 
\subsection{Graphs}

Let $K_n$ denote the complete graph with $n$ vertices, $P_n$ denote the $n$-vertex path, and $\MultiPartite{n_0,\,\dots,\,n_{m-1}}$ denote the complete $m$-partite graph with parts of sizes $n_0,\dots,n_{m-1}$. For integers $m$ and $n$, the \emph{Tur\'an graph} $\Turan{m}{n}$ is the complete $m$-partite graph on $[n]$ with almost equal parts, that is,   $\Turan{m}{n}\cong \MultiPartite{n_0,\,\dots,\,n_{m-1}}$ with $n_i=\lfloor (n+i)/m\rfloor$ for $i\in [m]$.
Also, $F\sqcup G$ denotes the union of vertex-disjoint copies of graphs $F$ and~$G$.

For a graph $G=(V(G),E(G))$, its \emph{order} is $v(G):=|V(G)|$.
The \emph{neighbourhood} of a vertex $u\in V(G)$ is the set 
\[
\Gamma_G(u):=\{w\in V(G)\mid \{u,w\}\in E(G)\}.
\]The \emph{degree} of $u$ is $\deg_G(u):=|\Gamma_G(u)|$.  Also, the \emph{complement} of a graph $G$ is
\[
\O G:=\left(V(G),\binom{V(G)}{2}\setminus E(G)\right)
\]
and, for $X\subseteq V(G)$, the subgraph \emph{induced} by $X$ is 
\[
G[X]:=(X,\{uw\in E(G)\mid u,w\in X\}).
\]
For distinct vertices $u,w\in V(G)$, let $G\oplus uw$ denote the graph obtained from $G$ by \emph{flipping} (i.e.\ changing) the adjacency between $u$ and $w$.

The \emph{edit distance} $\dedit(G,G')$ between two graphs $G$ and $G'$ of the same order is the minimum value of $|f(E(G))\bigtriangleup E(G')|$ over all bijections $f:V(G)\to V(G')$; in other words, it is the smallest number of \emph{edits} (changes in adjacency) we have to do in one graph to make it isomorphic to the other. The distance from a graph $G$ to a graph family $\C G$ is 
\[
\dedit(G,\C G):=\min\{\dedit(G,G')\mid G'\in\C G,\ v(G')=v(G)\}.
\]

For $m\in \I N$, let $\C F_m$ denote the family of graphs of order $m$ consisting of one representative from each isomorphism class. For graphs $F$ and $G$ with $\kappa$ and $n$ vertices respectively, let $\P(F,G)$ be the number of $\kappa$-subsets $X\subseteq V(G)$ that induce a subgraph isomorphic to $F$ in~$G$. 
If $\kappa\le n$ then we define $\p(F,G):=\P(F,G)/\binom{n}{\kappa}$, calling it the \emph{(induced) density} of $F$ in~$G$.

For a graph $B$ on $[m]$ and pairwise disjoint sets $V_0,\dots,V_{m-1}$ (some possibly empty), the \emph{blowup} $\blow{B}{V_0,\dots,V_{m-1}}$ of $B$ is the graph on $V=\cup_{i=0}^{m-1} V_i$ where $x\in V_i$ and $y\in V_j$ are adjacent if and only if $\{i,j\}\in E(B)$. In particular, each part $V_i$ is an independent set consisting of \emph{clones}, that is, (non-adjacent) vertices having the same sets of neighbours. Let $\blow{B}{}$ denote the family of all blowups of~$B$. We may call $B$ the \emph{base graph}.

\subsection{Extremal $\obj[\gamma]{}$-problem}\label{se:gamma}

Given an integer $\kappa\ge 2$ and a function
$\gamma:\C F_\kappa\to\I R$, we consider the following function on graphs:
\begin{equation}\label{eq:objGamma}
\obj[\gamma]{G}:=\sum_{F\in\C F_\kappa} \gamma(F)\p(F,G),\quad \mbox{for a graph $G$ with $v(G)\ge \kappa$}.
 \end{equation}
  The corresponding extremal \emph{$\obj[\gamma]{}$-problem} asks for the value of
 \begin{equation}\label{eq:objGammaN}
 \obj[\gamma]{n}:=\max\{\obj[\gamma]{G}\mid v(G)=n\},\quad \mbox{for $n\ge\kappa$}.
 \end{equation}
 Also, we define $\obj[\gamma]{}:=\lim_{n\to\infty}\obj[\gamma]{n}$, where the existence of the limit follows from a simple averaging argument, see for example~\PST{Lemma~2.2}.

The \emph{$(m-1)$-dimensional simplex} is
\beq{eq:Simplex}
 \I S_m:=\left\{(x_0,\dots,x_{m-1})\in \I R^m\mid x_0+\dots+x_{m-1}=1\mbox{ and } \forall\, i\in [m]\ x_i\ge 0\right\}.
 \eeq
For $\V x=(x_0,\dots,x_{m-1})\in \I S_m$ and a graph $B$ on $[m]$, let
$\obj[\gamma]{\blow{B}{\V x}}$ be the limit as $n\to\infty$ of $\obj[\gamma]{\blow{B}{V_0,\ldots,V_{m-1}}}$, where $|V_i|=(x_i+o(1))n$ for $i\in [m]$. This is a continuous function on $\I S_m$ (in fact, a polynomial). 
Let 
\[
\obj[\gamma]{\blow{B}{}}:=\sup\{\obj[\gamma]{\blow{B}{\V x}}\mid \V x\in\I S_m\}.
\]
By the compactness of $\I S_m$ the supremum is attained by at least one $\V x\in\I S_m$; such vectors will be called \emph{$(\obj[\gamma]{},B)$-optimal}.
The graph $B$ is  called \emph{$\obj[\gamma]{}$-minimal} if for every graph $B'$ obtained from $B$ by removing a vertex, it holds that $\obj[\gamma]{\blow{B'}{}}<\obj[\gamma]{\blow{B}{}}$. A standard compactness argument shows that $B$ is $\obj[\gamma]{}$-minimal if and only if every $(\obj[\gamma]{},B)$-optimal vector has all entries positive. Also, $B$ is called \emph{$\obj[\gamma]{}$-optimal} if $\obj[\gamma]{}=\obj[\gamma]{\blow{B}{}}$, that is, we can attain the asymptotically optimal constant $\obj[\gamma]{}$ by some blowups of~$B$.

Let us pause and see how some of these definitions relate to the semi-inducibility problem  for a 2-coloured graph $H$. Observe that 
the graph function $\obj[H]{}$ can be represented as $\obj[\gamma]{}$: take $\kappa:= v(H)$ and define the value of $\gamma$ on $F\in\C F_\kappa$ to be $\obj[H]{F}$ (which is the probability that a random bijection $V(H)\to V(F)$ gives an embedding of $H$ into~$F$). 
Whenever $H$ is given, these will be the default values of $\kappa$ and $\gamma$; then the functions $\obj[H]{\cdot}$ and $\obj[\gamma]{\cdot}$ will be the same.
Let us call a (not necessarily injective) map $f:V(H)\to V(B)$ a \emph{homomorphism}  from $H$ to $B$ and write it as $f:H\to B$ if the image of every red edge is an edge while every blue edge is sent to a non-edge or to two equal vertices of $B$; equivalently, it is a possible assignment of the vertices of $H$ to the parts under an embedding of $H$ into a blowup of $B$.
Using this notation, 
we have
 \beq{eq:ojbBX}
  \obj[H]{\blow{B}{\V x}}=\sum_{f:H\to B}\prod_{i\in V(H)} x_{f(i)},\quad \mbox{for $\V x\in\I S_m$}.
  \eeq

We will also need some definitions  from \PST{}.

A graph $B$ on $[m]$ is \emph{$\obj[\gamma]{}$-flip-averse} if there is $\delta>0$ such that, for every blow-up $G=\blow{B}{V_0,\dots,V_{m-1}}$ with $n\ge 1/\delta$ vertices and $\obj[\gamma]{G}\ge \obj[\gamma]{\blow{B}{}}-\delta$, it holds for any distinct vertices $u,w\in V(G)$ (possibly from the same part) that
 \begin{equation}\label{eq:flip}
 \obj[\gamma]{G}-\obj[\gamma]{G\oplus uw}\ge \delta n^{-2}.  
 \end{equation}
Informally, this definition states that if we flip any pair in a large (almost) optimal blowup of $B$ then we decrease the objective function by the maximum possible amount in the order of magnitude. By compactness of $\I S_m$, this can be equivalently re-stated in terms of the optimal vectors in~$\I S_m$. For example, $B$ is $\obj[H]{}$-flip-averse if and only if, for every $(\obj[H]{},B)$-optimal vector $\V x\in \I S_m$ and for every (not necessarily distinct) vertices $u,w\in V(B)$, it holds that
 \begin{equation}\label{eq:flipX}
 \sum_{i\in V(H)} \sum_{j\in V(H)\setminus\{i\}} \left(\sum_{f:H\to B\atop f(i)=u,\ f(j)=w} \prod_{\ell\in V(H)\setminus\{i,j\}} x_{f(\ell)}
 -\sum_{f:H\to B\oplus uw\atop f(i)=u,\ f(j)=w} \prod_{\ell\in V(H)\setminus\{i,j\} }
 x_{f(\ell)}\right)>0.
 \end{equation}
 If we know the set of optimal vectors and it is finite, then~\eqref{eq:flipX} can be checked on computer.
 
Next, we give the definition, also from~\PST{}, when an extremal problem is strict. Roughly speaking, it states that if we take a large optimal blowup $G$ of $B$ and attach a new vertex $u$ then the normalised contribution of $u$ to the objective function is equal to $\obj[\gamma]{}+o(1)$ (the value observed for vertices in $G$) only if $u$ is almost a clone of an existing vertex. Formally, 
a graph $B$ on $[m]$ is called \emph{$\obj[\gamma]{}$-strict} if $B$ is $\obj[\gamma]{}$-minimal and, for every $(\obj[\gamma]{},B)$-optimal $\V x\in\I S_m$ and every vector $\V y\in [0,1]^m$, the polynomial $p(\V x,\V y)$ defined below in~\eqref{eq:StrictP} assumes value $\obj[\gamma]{}$ only if there is $i\in [m]$ with each $y_j$ being the indicator function of $ij\in E(B)$. To define the polynomial $p$, let $B'$ be the blow-up of $B$ with parts $V_i:=\{i,i+m\}$ for $i\in [m]$ and let $B''$ be obtained from $B'$ by adding $2m$ as a new vertex with $\Gamma_{B''}(2m)=[m]$. 
(In other words, we add a clone of each vertex of $B$ and then add a new vertex adjacent to the original vertices of~$B$.)  We define 
\begin{equation}\label{eq:StrictP}
 p(\V x,\V y):=\sum_{F\in\C F_\kappa} \gamma(F) \sum_{f}\prod_{j\in V(F)\setminus f^{-1}(2m)} z_{f(j)},
\end{equation}
  where 
  \[
  \V z:=(y_0x_0,\dots,y_{m-1}x_{m-1},(1-y_0)x_0,\dots,(1-y_{m-1})x_{m-1})
  \]and sum is taken over all (including non-injective) maps $f:V(F)\to V(B'')$ that send edges to edges and non-edges to non-edges while the pre-image of  the special vertex $2m$ consists of exactly one vertex.
  
In the special case when $\obj[\gamma]{}=\obj[H]{}$, we have that \[
 p(\V x,\V y)=\sum_{f:H\to B''\atop |f^{-1}(2m)|=1}\prod_{j\in V(H)\atop f(j)\not=2m} z_{f(j)}.
  \]
 We refer the reader to Claim~\ref{cl:H4strict} to see an example where the strictness property is verified.
 
Checking strictness can be performed by a computer calculation, if we already have all $(\obj[\gamma]{},B)$-optimal vectors $\V x$, as follows. Assume that none of these vectors has a zero entry (as otherwise the graph $B$ is not $\obj[\gamma]{}$-minimal and thus is not $\obj[\gamma]{}$-strict by definition). For each $(\obj[\gamma]{},B)$-optimal vector $\V x$, we have to find all elements $\V y\in [0,1]^m$ for which the polynomial $p(\V x,\V y)$ attains its maximum value~$\obj[\gamma]{}$. We run over all $3^n$ choices of assigning to each $y_i$ value $0$, value $1$ or leaving $y_i$ free. For each choice, we use Buchberger's algorithm to check if there are vectors $\V y$ as above such that $p(\V x,\V y)=\obj[\gamma]{}$, $\partial p(\V x,\V y)/\partial y_i=0$  for each free variable $y_i$, $\V y\in [0,1]^m$ but $\V y$ is not the $\{0,1\}$-valued vector encoding  the $B$-adjacency of some $i\in [m]$.

For $\V x\in\I S_m$, we call the $\obj[\gamma]{}$-problem  \emph{Erd\H os--Simonovits $(B,\V x)$-stable} if for every $\e>0$ there are $\delta>0$ and $n_0$ such that if $G$ is a graph with $n\ge n_0$ vertices and $\obj[\gamma]{G}\ge \obj[\gamma]{}-\delta$ then $G$ is $\e {n\choose 2}$-close in the edit distance to a blowup $\blow{B}{V_0,\dots,V_{m-1}}$ with $\left|\,|V_i|-x_in\,\right|<\e n$. Recall that the last statement means that there is a partition $V(G)=V_0\cup\dots\cup V_{m-1}$ such that
 \beq{eq:ESStab}
 \left|{\textstyle E(G)\bigtriangleup} E(\blow{B}{V_0,\dots,V_{m-1}})\right|\le \e {n\choose 2}.
 \eeq

This property is very useful as the first step towards characterising graphs of sufficiently large order $n$ with $\obj[\gamma]{G}=\obj[\gamma]{n}$. This approach was pioneered by Erd\H os~\cite{Erdos67} and
Simonovits~\cite{Simonovits68}.

Following~\PST{},  we call the $\obj[\gamma]{}$-problem \emph{perfectly $B$-stable}\label{def:PerfectStability}
for a graph $B$ if there is $C>0$ such that for every graph $G$ of order $n\ge C$ there is a blowup $G'$ of $B$ such that 
 \beq{eq:PerfectStabDef}
 \dedit(G,G')\le C\,
 \left(\obj[\gamma]{n}-\obj[\gamma]{G}\right)n^2.
 \eeq
In particular, it follows that, for all $n\ge C$, every order-$n$ graph $G$ with $\obj[\gamma]{G}=\obj[\gamma]{n}$ is a blowup of $B$. 
Perfect stability has been established for various extremal problems, see e.g.~\cite{Furedi15,NorinYepremyan17,NorinYepremyan18,RobertsScott18,PikhurkoSliacanTyros19,LiuPikhurkoSharifzadehStaden23}. 

Let us remark that the property of being perfectly $B$-stable remains true if we add new vertices to $B$ with arbitrary attachments. If $B^+$ is obtained from $B$ in this way, then setting the new part sizes to be empty, we get that $\blow{B}{V_0,\dots,V_{m-1}}$ is isomorphic to
$\blow{B^+}{V_0,\dots,V_{m-1},\emptyset,\dots,\emptyset}$.
In all cases of this paper when we establish perfect stability, the listed graph $B$ is $\obj[\gamma]{}$-minimal, that is, has no redundant vertices.

Since $\gamma$ or $H$ will be fixed and clear from the context, we may omit them from our notation. Asymptotic notation, such as $o(1)$, is taken with respect to $n\to\infty$ (where $n$ is usually the order of the unknown graph $G$); the constants hidden in it may depend on $\kappa$ and $\gamma$ but not on any other parameters. We call a sequence of graphs $(G_n)_{n\in\I N}$ with increasing orders \emph{almost extremal} (resp.\ \emph{almost $c$-regular}) if $\obj[\gamma]{G_n}=\obj[\gamma]{}+o(1)$ as $n\to\infty$ (resp.\ for every $\e>0$ there is $n_0$ such that, for every $n\ge n_0$, at least $(1-\e)v(G)$ vertices $u$ of $G$ satisfy $|\deg_G(u)/v(G)-c|\le \e$).

\subsection{Flag Algebras}

Since the flag algebra approach of Razborov~\cite{Razborov07}
is well established by now (and is described in detail in e.g.~\cite{Razborov10,BaberTalbot11,SFS16,GilboaGlebovHefetzLinialMorgenstein22}), we give only high-level description. Suppose we are given an integer $\kappa$ and a function $\gamma:\C F_{\kappa}\to\I R$ and would like to upper bound the corresponding extremal function $\obj[\gamma]{n}$ defined in~\eqref{eq:objGammaN}. A flag algebra proof of $\obj[\gamma]{n}\le u+o(1)$ amounts to an identity that holds for every graph $G$ of order $n\to\infty$ of the form
 \beq{eq:FAMain}
  u-\sum_{F\in\C F_\N} \obj[\gamma]{F}\p(F,G)=\mathrm{SOS}+ \sum_{F\in\C F_\N} c_F\,\p(F,G)+o(1), 
  \eeq
where $\N\ge \kappa$ is an integer,  each $c_F$ is a non-negative real, called the \emph{slack} of $F\in \C F_\N$ while $\mathrm{SOS}$ is a sum of squares (of certain rooted densities) which can also be equivalently written as $\sum_{F\in\C F_\N} a_F\p(F,G)+ o(1)$ for some reals~$a_F$. 
(See e.g.~\eqref{eq:H11SOS} for an example of an $\mathrm{SOS}$-expression in the flag algebra context.) 
The corresponding \emph{certificate} lists $u$, $\N$, the slacks $c_F$ and the positive semi-definite matrices that define the  $\mathrm{SOS}$-term. All of these are independent of $n$, and are stored using exact arithmetic.

All presented certificates were generated using the package \texttt{FlagAlgebraToolbox} by the first author \cite{flagalgebratoolbox}. The scripts that we use to generate the certificates, verify the identity in~\eqref{eq:FAMain} as well as any additional claimed properties can be found 
in a separate GitHub repository \href{https://github.com/bodnalev/supplementary_files}{\url{https://github.com/bodnalev/supplementary_files}} inside the folder \verb$graph_semi_inducibility$, while the certificates themselves can be found in its sub-folder \verb$certificates$. Alternatively, all these files can be found in the ancillary folder of the arXiv version of this paper.
Each certificate starts with \verb$semiind$ plus the index as in Table~\ref{ta:semiind}.
Alternatively, a reader can use their own verifier; a description of how the data are arranged in each certificate file can be found in the \href{https://github.com/bodnalev/supplementary_files/blob/main/readme.md}{readme.md} file.

\section{Individual cases}\label{se:proofs}

Here we investigate individual graphs $H_i$ from Table~\ref{ta:semiind}. First, we list those cases (namely $H_0$, $H_1$, $H_4$ and $H_{15}$) where we were able to prove all results without using flag algebras. Our proof for the graph $H_5$ requires some basics of the limit theory of complete partite graphs, so we postpone it to the end.

 \newcommand{\lambdamin}{\lambda_{\mathrm{min}}}
 \newcommand{\lambdamax}{\lambda_{\mathrm{max}}}
 \newcommand{\gammamax}{\gamma_{\mathrm{max}}}
 \newcommand{\supp}{\mathrm{supp}}
 \newcommand{\eps}{\varepsilon}
 \newcommand{\OPT}{\mathrm{OPT}}
 \newcommand{\Syone}{(Sym1)}
 \newcommand{\Sytwo}{(Sym2)}
 \newcommand{\Sone}{(Str1)}
 \newcommand{\Stwo}{(Str2)}
 
 \newcommand{\Dedit}{\Delta_{\mathrm{edit}}}
 \newcommand{\Done}{\hat{\Delta}_{1}}
 \newcommand{\done}{\hat{\delta}_{1}}
 \newcommand{\PC}{\overline{\mathcal{P}}}

\subsection{Stars and matchings}
 
Here we provide two easy observations (that in particular apply to the graphs $H_0$ and $H_1$ of Table~\ref{ta:semiind}). 
They follow from the facts that $x_0:=\ell/(\ell+m)$ is the unique maximiser of $x^\ell(1-x)^m$ on the interval $[0,1]$ and that, since we maximise a continuous function on a compact set, all $x\in [0,1]$ that almost attain the maximum value must be close to~$x_0$.

 \begin{observation}\label{ob:SIMatching} If $H$ is the matching of size $\ell+m$ with $\ell$ red and $m$ blue edges then 
  \[
   \obj[H]{}= \frac{\ell^\ell m^m}{(\ell+m)^{\ell+m}},
   \]
   and a large graph is almost extremal if and only if its edge density is $\ell/(\ell+m)+o(1)$.\qed\end{observation}

 \begin{observation}\label{ob:SIStar} If $H$ is the star $\MultiPartite{\ell+m,\,1}$ with $\ell$ red and $m$ blue edges then 
  \[
   \obj[H]{}= \frac{\ell^\ell m^m}{(\ell+m)^{\ell+m}},
   \]
   and a large graph is almost extremal if and only if it is almost $\ell/(\ell+m)$-regular.\qed\end{observation}

\subsection{Exact result for $H_0$}\label{se:H0}

Observation~\ref{ob:SIMatching} gives that $\obj[H_0]{}=1/4$ and describes all almost extremal graphs. In this section, we determine the extremal function exactly. For integers $n$ and $0\le e\le \binom{n}2$, let the \emph{quasi-clique} $QC(n,e)$ be obtained from the clique $K_{k+1}$ by removing $\binom{k+1}2-e$ edges forming a star and then adding $n-k-1$ isolated vertices, where $k$ is the maximum integer with $\binom k2< e$. Clearly, $QC(n,e)$ has $n$ vertices and $e$ edges. 

\begin{theorem}\label{th:H0exact} Let \[H := H_0 = \DrawSemiGraphLabeled{4}{0/1}{2/3}{0,2,1,3},\] which is the graph on $[4]$ with a red edge $02$ and a blue edge $13$. Take any $n\ge 4$. Define $k_0:=\lfloor x_0\rfloor$, where
\begin{equation}\label{eq:x0}
 x_0:=\frac{3+\sqrt{2n^2-10n+13}}2.
 \end{equation}
 Then $\Obj{n}=\Obj{QC(n,\binom{k_0+1}2)}$. Moreover, if  $2n^2-10n+13$ is not a square then $QC(n,\binom{k_0+1}2)$ and its complement are the only extremal graphs; otherwise there are at least $2(k_0+1)$ non-isomorphic extremal graphs.
 \end{theorem}
 
\bpf The easy case $n=4$ (when $\Obj{n}=12$ and there are exactly two extremal graphs) satisfies the theorem, so assume  that $n\ge 5$.

First, let us prove the stated upper bound on $\Obj{n}$. Take any extremal graph $G$ with $n$ vertices. Let $e$ denote the number of edges in~$G$. We have 
\begin{equation}\label{eq:DegSquared}
\Obj{G}=4e\left(\binom{n}2-e\right)-4\sum_{u\in V(G)} \deg(u)(n-1-\deg(u)),
\end{equation}
since every pair consisting of an edge $uv$ and a non-edge $wx$ gives 4 embeddings of $H$, except we have to exclude intersecting pairs (i.e.\ those with $\{u,v\}\cap \{w,x\}\not=\emptyset$). Thus, if we fix the number of edges of $G$ then our objective is to maximise the sum of the squares of degrees. Ahlswede and Katona~\cite[Theorem~2]{AhlswedeKatona78} proved that, for every $(n,e)$, the maximum value is attained by $QC(n,e)$ or by the complement of $QC(n,\binom{n}{2}-e)$. Since $\Obj{G}=\Obj{\O G}$, let us assume that $G=QC(n,e)$. Let $k$ be the maximum integer with $\binom k2< e$ and let $\ell:=e-\binom{k}2$. Thus, apart from isolated vertices, $G$ can be viewed as obtained from the clique on a $k$-set $A$ by adding $\ell$ edges, $1\le \ell\le k$, all incident to the same vertex~$u$ outside of $A$. If we fix $n$ and $k$  then $\Obj{G}$ is a linear function of $\ell$ with the coefficient at $\ell$ being
 \begin{equation}\label{eq:Cnk}
 c_{n,k}:=\binom{n-k-1}2+(n-k-1)(k-1)-\binom{k-1}2=\frac{n^2-5n}2-k^2+3k+1.
 \end{equation}
 We can assume that $\ell=k$, for otherwise $c_{n,k}=0$ (since we can both increase and decrease $\ell$) and thus $\Obj{G}$ does not depend on~$\ell$. Thus $G=K_{k+1}+\O K_{n-k-1}$. We have that \[
 \Obj{G}-\Obj{K_k+\O K_{n-k}}=k\,c_{n,k}.
 \] 
 The function $c_{n,k}$, for a real variable $k$, has a root at $k=x_0$, is negative for $k>x_0$ and positive for $0<k<x_0$ (since $n\ge 5$). Thus $\Obj{n}=\Obj{G}\le \Obj{QC(n,\binom{k_0+1}2)}$, proving the first claim of the theorem.
 
Now, let us investigate the structure of extremal graphs. Take any extremal graph $G$ of order~$n$. By~\eqref{eq:DegSquared}, it also maximises $S(G):=\sum_{u\in V(G)} \deg^2(u)$ given its size~$e=e(G)$. 
It is easy to see that $G$ has either an isolated or a universal vertex. Indeed, if $u$ is a vertex of the minimal degree and $w$ is a vertex of the maximal degree and, say, $uw\in E(G)$ then $\deg(w)=n-1$: if there existed a non-neighbour $v$ of $w$ then we could strictly increase $S(G)$ by replacing the edge $uw$ by~$vw$. Using a version of this argument, it can be derived via a suitable induction statement that $G$ is a \emph{threshold graph}, that is, it can be obtained from a single vertex by adding new isolated and universal vertices in some order. 
While the result of Ahlswede and Katona~\cite{AhlswedeKatona78} shows that at least one graph that maximises $S(G)$ for given $(n,e)$  has very simple structure, the task of describing all extremal graphs for this problem turns out to be rather complicated, even though it is enough to consider threshold graphs only. This task is addressed in \cite{Byer99,PeledPetreschiSterbini99,AbregoFernandezNuebauerWatkins09}. In fact, there are six families that between them contain every extremal graph.

Let us return back to our problem. By taking complements, we can assume that $G$ has an isolated vertex.
Note that $c_{n,k}$ as a function of real $k\ge0$ is zero only if $k=x_0$, where $x_0$ and $c_{n,k}$ are defined in~\eqref{eq:x0} and~\eqref{eq:Cnk} respectively.

Suppose first that $2n^2-10n+13$ is not a square. Then $c_{n,k}$ is non-zero for every integer~$k\ge 0$. Thus the parameter $\ell$ in the proof of the upper bound is necessarily $k$, that is, $e=\binom{k_0+1}2$. For such a special value of $e$ (plus our assumption that $G$ has an isolated vertex), Theorem~2.4 in \cite{AbregoFernandezNuebauerWatkins09} implies that $G$ is isomorphic to $QC(n,\binom{k_0+1}2)$, as desired.

Finally, suppose that $2n^2-10n+13$ is a square. Then $x_0=k_0$ is an integer and it holds that $k_0\le n-2$ (by $n\ge 6$) and $c_{n,k_0}=0$. Thus  $QC(n,\binom{k_0}2+\ell)$ for $0\le \ell\le k_0$ are all pairwise non-isomorphic extremal graphs (which, by containing an isolated vertex, are non-isomorphic to any of their complements).\epf

Let us remark that there are infinitely many $n$ such that $2n^2-10n+13$ is a square of some integer $d$. The corresponding Diophantine equation can be re-written as $2d^2-1=m^2$, using the substitution $m=2n-5$. This a special case of the well-studied negative Pell equation and an explicit description of all (infinitely many) pairs of integers $(d,m)$ with $2d^2-1=m^2$ can be found in e.g.\
\cite[Page 144]{AndreescuAndricaCucurezeanu10ide}. For such $n$ (when we are free to choose any $0\le \ell\le k_0$), there are choices of $\ell$ when some other extremal graphs appear; we refer the reader to Cases~2.2 and 2.3 of \cite[Theorem~2.4]{AbregoFernandezNuebauerWatkins09} for details.

\subsection{Exact result for $H_1$}\label{se:H1}
 
For $H_1$, the exact function can be easily determined. Recall that $\pi_k(n)$, defined in~\eqref{eq:pi}, is the maximum product of $k$ non-negative integers summing up to~$n$.

\begin{lemma} Let \[H := H_1 = \DrawSemiGraphLabeled{4}{0/1}{0/2,0/3}{0,1,2,3},\] which is the graph on $[4]$ that consists of a red edge $01$ and two blue edges $02$ and $03$. Take any integer $n\ge 4$. 

If $n$ is congruent to $1$, $2$ or $4$ modulo 6 then it holds that $\Obj[H_1]{n}= n\,\pi_3(2n-3)/2$ and a  graph $G$ of order $n$ attains equality if and only if $G$ is $\lfloor n/3\rfloor$-regular.

If $n=6k$ for an integer $k$ then $\Obj[H_1]{n}= n\,\pi_3(2n-3)/2-3k(4k-1)$ and a  graph $G$ of order $n$ attains equality if and only if every vertex has degree $2k-1$ or $2k$.

If $n=6k+3$ for an integer $k$ then $\Obj[H_1]{n}= n\,\pi_3(2n-3)/2-(6k+3)(4k+1)/2$ and a  graph $G$ of order $n$ attains equality if and only if every vertex has degree $2k$ or $2k+1$.

Finally, if $n=6k+5$ for an integer $k$ then $\Obj[H_1]{n}= n\,\pi_3(2n-3)/2-4k-2$ and a  graph $G$ of order $n$ attains equality if and only it has $n-1$ vertices of degree $2k+1$ and 1 vertex of degree $2k+2$.
 \end{lemma}
 
 \bpf The $H_1$ count in an $n$ vertex graph $G$ can be calculated by summing over the choices $u_0 \in V(G)$ mapping to the center of the star, the choices of $u_1 \in \Gamma(u_0)$ and distinct $u_2, u_3 \in \Gamma_{\O G}(u_0)$, which can be done in $$\deg(u_0) \big( n-1-\deg(u_0) \big) \big( n-2-\deg(u_0)  \big)$$ different ways. 
 Given an integer $n\ge 1$, let $r:=\lfloor n/3\rfloor$ and define 
 \[
 D(d):=\pi_3(2n-3)-2d\,(n-1-d)(n-2-d),\quad \mbox{for $d\in\I N$}.
 \]
 Clearly, it holds that $D(d)\ge 0$ for any integer $d\in[n]$. Thus, for an $n$-vertex graph $G$, we have that 
 \begin{equation}\label{eq:H1}
  \Obj[H_1]{G}
  =\frac12 \sum_{u_0\in V(G)}\left(\pi_3(2n-3)-D(\deg(u_0))\right)\le \frac{n\,\pi_3(2n-3)}{2}.
  \end{equation}

If $n$ is $1$,  $2$ or $4$ modulo 6 then, as it is routine to see, $G$ attains equality in~\eqref{eq:H1} if and only if it is $r$-regular (and at least one such graph exists since $nr$ is even is each of these cases); thus the lemma holds.

Next, suppose that $n=6k$. Here $2n-3=12k-3$. We have that $\pi_3(12k-3)=(4k-1)^3$. Furthermore, the product of three non-negative integers, not all equal, summing up to $12k-3$ is at most $(4k-2)(4k-1)4k=\pi_3(12k-3)-4k+1$, with equality if and only if these numbers are $4k-2$, $4k-1$ and $4k$ in some order. Thus, for $d\in [n]$, we have that $D(d)\ge 4k-1$ with equality if and only if $2d$ makes a triple of consecutive integers with $6k-2-d$ and $6k-1-d$, that is, 
$d=2k-1$ or $d=2k$. Thus $\sum_{u_0\in V(G)}D(\deg(u_0))\ge n\cdot (4k-1)$, giving the stated upper bound on $\obj[H_1]{G}$, with equality if and only if every degree is $2k-1$ or $2k$. Furthermore, the set of such graphs is non-empty (e.g.\ by the Erd\H os--Gallai theorem~\cite{ErdosGallai60}).

The case $n=6k+3$ is completely analogous to the previous one, so we move to the remaining case $n=6k+5$ with $k\in\I N$. Here $2n-3=12k+7$ and $\pi_3(12k+7)=(4k+2)^2(4k+3)$. For $d\in [n]$, the product $2d\,(6k+3-d)(6k+4-d)$ is equal to $\pi_3(12k+7)$ if and only if $2d=6k+3-d$, that is, $d=2k+1$. However, $(2k+1)n$ is odd so at least one vertex of $G$ has degree different from $2k+1$. As it is routine to see, the minimum of $D(d)$ for $d\in [n]\setminus \{2k+1\}$ is $8k+4$ and it is attained if and only if $d=2k+2$. All claims of the lemma in this case routinely follow from the above observations.\epf

\subsection{Graph $H_4$}\label{se:H4}

Here we resolve the case of $H_4$, with perfect stability.

\begin{theorem} Let \[H := H_4 = \DrawSemiGraphLabeled{4}{0/2,0/3,1/2}{0/1}{0,1,3,2},\] that is, the graph on $[4]$ with the red edge set $\{02,03,13\}$ and the blue edge set $\{01\}$. Then all of the following statements hold.
\begin{enumerate}
\item\label{it:H4obj} $\obj{}=4/27$.
\item\label{it:H4stab} The semi-inducibility problem for $H$ is Erd\H os-Simonovits $(K_3,(1/3,1/3,1/3))$-stable.
\item\label{it:H4perfect} The semi-inducibility problem for $H$ is perfectly $K_3$-stable.
\end{enumerate}
\end{theorem}
\bpf  The 3-partite  Tur\'an graphs $\Turan{3}{n}$ with $n\to\infty$ show that $\obj{}\ge 4/27$. To prove the matching upper bound, take an arbitrary graph $G$ with $n\ge 4$ vertices. By considering where $0\in V(H)$ is mapped, we derive that
 \begin{equation}\label{eq:H4Upper}
 \Obj{G}\le \sum_{u_0\in V(G)} (\deg(u_0))^{\fallingfact{2}} \cdot(n-1-\deg(u_0))\le \sum_{u_0\in V(G)} \frac{\pi_3(2n-2)}2 = \frac{n\,\pi_3(2n-2)}2.
 \end{equation}
 This is $(4/27+o(1))n^4$ as $n\to\infty$, proving the desired upper bound on $\obj{}$. 
 
As a side remark, observe that if we let $G$ be the Tur\'an graph $\Turan{3}{n}$ then all inequalities in~\eqref{eq:H4Upper} are equalities if $n$ is divisible by 3 (and we have a gap of only $O(n^2)$ for general~$n$). 

Let us turn to proving the Erd\H os-Simonovits stability. Let $n\to\infty$ and $G$ be any graph with $\obj{G}=(4/27+o(1))n^4$. Recall that $\O{P}_3=([3],\{01\})$ is the graph with $3$ vertices and exactly one edge.

\begin{claim}\label{cl:H4OP3}
 $p(\O{P}_3,G)=o(1)$.\end{claim}
 
\bpf[Proof of Claim.] Suppose that the claim is false. Then there are $\Omega(n)$ vertices $u_0$ for which the number of copies of $\O{P}_3$ in $G$ with $0$ sent to $u_0$ is $\Omega(n^3)$. Each such vertex $u_0$ has degree $\Omega(n)$ and thus contributes a gap of $\Omega(n^2)$ to~\eqref{eq:H4Upper}. Indeed, when counting embeddings $f$ of $H$ into $G$ with $f(0)=u_0$, we upper bound the number of choices of the pair $(f(1),f(3))$ by $\deg(u_0)(n-1-\deg(u_0))$; however $\Omega(n^2)$ of such pairs (those that create a copy of $\O{P}_3$ with $u_0$) are non-adjacent and do not extend to an embedding of $H_4$. Thus the total gap in~\eqref{eq:H4Upper} is $\Omega(n^4)$, a contradiction.\epf
 
By the Induced Removal Lemma~\cite{AlonFischerKrivelevichSzegedy00}, we can change $o(n^2)$ adjacencies in $G$ so that no copy of $\O{P}_3$ occurs, resulting in a complete multi-partite graph. 

\begin{claim}\label{cl:3Parts}
After removing $o(n)$ vertices, $G$ has three parts, each of size $(1/3+o(1))n$.
\end{claim}
\bpf[Proof of Claim.] Let the parts of $G$ be $V_0,\dots,V_{m-1}$ of sizes $x_0n\ge \ldots\ge x_{m-1}n$. Let $\ell$ be the minimal index such that $x_{\ell}<2/7$ (with $\ell:=m$ if $x_{m-1}\ge 2/7$). Clearly, $0\le \ell \le \lfloor 7/2\rfloor=3$. The argument leading to the upper bound~\eqref{eq:H4Upper} gives that
 \[
 \obj{G}\le \sum_{i=0}^{m-1} x_i\cdot x_i(1-x_i)^2+o(1)\le \frac{4}{27}\,\sum_{i=0}^{\ell-1} x_i + \frac27\left(1-\frac{2}{7}\right)^2\cdot  \sum_{i=\ell}^{m-1} x_i+o(1).
 \]
 By $\obj{G}=4/27+o(1)$ and $\sum_{i=0}^{m-1}x_i=1$, this implies that $\sum_{i=\ell}^{m-1} x_i=o(1)$ and that, for each $0\le i<\ell$, we have $x_i=1/3+o(1)$. 
 It follows that $\ell=3$. Thus we can make the graph $G$ 3-partite with each part of size $(1/3+o(1))n$, by removing all vertices in $\cup_{i=\ell}^{m-1}V_i$, as desired.\epf

Thus we can make $G$ into $\Turan{3}{n}$ by editing $o(n^2)$ adjacencies, which is exactly the claimed Erd\H os-Simonovits stability. 

\begin{claim}\label{cl:H4strict} The $\obj[H]{}$-problem is  $K_3$-flip-averse and $K_3$-strict.\end{claim}
\bpf[Proof of Claim.] By Claim~\ref{cl:3Parts}, the unique optimal vector for $K_3$ is $(1/3,1/3,1/3)$. Thus the flip averseness follows by computing that, for $n\to\infty$,
\[
\obj{\Turan{3}{n}}-\obj{\Turan{3}{n}\oplus uw}+o(1)=\left\{
\begin{array}{ll} \frac49+\frac49+\frac49-\frac49,& \mbox{if $uw\in E(G)$,}\\
-\frac49-0-0+\frac89,& \mbox{if $uw\in E(\O G)$,}
\end{array}\right. >0,
 \]
 where the four terms in each row correspond to the total contribution of maps $V(H)\to V(\Turan{3}{n})$ such that $uw$ is the image of $02$, $03$, $13$ and $01$ respectively. Flipping $uw$ can affect only those embeddings for which the pair $\{u,w\}$ is the image of one of the coloured edges of $H$. If $uw$ is an edge of $\Turan{3}{n}$, then flipping it destroys the embeddings in which $uw$ is the image of one of $02$, $03$ or $13$, and creates embeddings in which $uw$ is the image of the blue edge $01$, justifying the signs in the first row. If $uw$ is a non-edge, an analogous reasoning applies. 

In order to check strictness, we write the polynomial $q(\V y):=p((1/3,1/3,1/3),\V y)$  as the sum $\sum_{i=0}^3 q_i(\V y)$, where $p$ is as in~\eqref{eq:StrictP} and, for $i\in [4]$, $q_i$ corresponds to those functions $f$ that map $i$ to the special vertex. Thus $q_i$ is the limit as $n\to\infty$ of $1/n^{\fallingfact{3}}$ times the number of embeddings of $H$ into $G'$ that send $i\in V(H)$ to $u$, where $G'$ is obtained from the Tur\'an graph
    $\Turan{3}{n}$ by adding a vertex $u$ with $y_in/3$ neighbours in the $i$-th part. It is easy to see that
\begin{eqnarray*}
    q_0(\V y)& =& \frac{y_0+y_1+y_2}3\cdot\sum_{i\in [3]}\,\sum_{j\in [3]\setminus\{i\}} \frac{y_i}3\cdot\frac{1-y_j}3,\\
   q_1(\V y) &=& \sum_{i\in [3]}\,\sum_{j\in [3]\setminus\{i\}} \frac{y_i}3\cdot \frac{1-y_j}3\cdot \frac23,\\
    q_2(\V y) &=& \frac{y_0+y_1+y_2}3\cdot \frac23\cdot \frac13,\\
    q_3(\V y)&=&\frac{y_0^2+y_1^2+y_2^2}9\cdot \frac23.
\end{eqnarray*}
Adding them, we obtain the (symmetric) polynomial $q(\V y)$. It can be checked to satisfy $q(0,1,1)=4\obj[H]{}=16/27$ (as expected) and our task is to show this element and its  permutations are the only elements of $[0,1]^3$ that attain the maximum value $16/27$.
Routine calculations show that, for any $i,j\in[3]$, it holds that
 \[
 \frac{\partial q(\V y)}{\partial y_i}-\frac{\partial q(\V y)}{\partial y_j}=
 \frac2{27}\cdot (y_i-y_j)\cdot (4-y_0-y_1-y_2). 
 \]
 This is zero for $\V y\in [0,1]^2$ only if $y_i=y_j$.
 Thus the only points that we need to check are those for which each coordinate is 0, 1 or some common value $y\in (0,1)$. There are a few cases to consider. For each case, we have to find the values of a univariate polynomial of degree at most 3 on its critical points in $(0,1)$ and check that we obtain $16/27$ only when the corresponding full vector $\B y\in [0,1]^3$ is a permutation of $(0,1,1)$. 
 This is a routine task, so we just do one case for an illustration. 
 For example, one can compute that $q(y,y,y)=(-6 y^3 + 4 y^2 + 6 y)/9$ and its derivative has only one root in $(0,1)$, namely $y=(2+\sqrt{31})/9$ with the value of the polynomial being $(356+124\sqrt{31})/2187<16/27$, so we do not get any counterexample to strictness in this case. We refer the reader to our code that does full verification as described in Section~\ref{se:Prelim} (without using any symmetries of the base graph).\epf

Theorem~5.8 in~\PST{} states, modulo minor technical assumptions that hold here, that if, with respect to some base graph $B$, an extremal $\obj[\gamma]{}$-problem is Erd\H os-Simonovits $B$-stable while the graph $B$ is both $\obj[\gamma]{}$-flip-averse and $\obj[\gamma]{}$-strict then the problem is perfectly stable. It applies to the $H_{4}$-semi-inducibility by the above claims, thus proving perfect $K_3$-stability.\epf

\subsection{Graph~$H_{15}$}
\label{se:H15}
  
  Recall that $H_{15}$ has red edges $01,02$ and blue edges are $03,12,13$.    
 Here we determine the exact value of $\Obj[H_{15}]{n}$ except when $n$ is 1 or 4 modulo 6 when we will have a gap of order $n^2$ between our bounds. 
  
As a building block we will need a construction of a triangle-free $d$-regular $n$-vertex graph for some special pairs $(n,d)$. The following two simple constructions will be our basic blocks. For integers $m\ge d\ge 0$, let $B_{m,d}$ be the bipartite graph with $2m$ vertices where each part is a copy of the residues modulo $m$ and we connect two residues across if their sum modulo $m$ is in $[d]$. Clearly, the graph $B_{m,d}$ is triangle-free (as it is bipartite) and $d$-regular. Next, for any positive integers $m,k$ with  $3k< m+2$, let $R_{m,k}$ be the graph whose vertex set consists of the residues modulo $2m+1$ and two vertices are adjacent if their difference modulo $2m+1$ belongs to the (symmetric under negation) interval 
  \[
  I:=\{m+1-k,m+2-k,\dots,m-1+k,m+k\}.
  \] 
   In other words, if we view the vertices as arranged on a circle then we connect each vertex to the diametrally opposite interval of size~$2k$.  Clearly, this graph is $2k$-regular.
   In order to check that there are no triangles, it is enough, by the symmetries of the construction, to check that two cyclically shortest hops from 0 jump over the interval $I$, that is, we have $2(m+1-k)>m+k$, which is exactly the assumption that we made.

\begin{theorem}\label{th:H15Exact} Let \[H := H_{15} = \DrawSemiGraphLabeled{4}{0/2,0/3}{0/1,1/2,2/3}{0,3,1,2},\] which is the graph on $[4]$ with red edge set $\{01, 02\}$ and blue edge set $\{03, 12, 13\}$. 
For every $n\ge 4$, it holds that $\Obj{n}\le n\,\pi_3(n-1)$.

If $n$ is 0, 2, or 5 modulo 6 then $\Obj{n}= n\,\pi_3(n-1)$ and an order-$n$ graph $G$ is extremal  if and only if $G$ is $\lceil n/3\rceil$-regular and triangle-free.

  If $n=6k+3$  for $k\in\I N$ then $\Obj{n}=n\,\pi_3(n-1)-2k$ and an order-$n$ graph $G$ is extremal if and only if $G$ is triangle-free, has exactly $n-1$ vertices of degree $2k+1$, and one vertex of degree $2k+2$.

If $n=6k+1$ or $n=6k+4$ for integer $k\to\infty$ then
\begin{equation*}6k^2+O(k)\le n\,\pi_3(n-1)-\Obj{n}\le \frac{15}2\,k^2+O(k).
\end{equation*}

In particular, it follows from above that $\obj{}=1/27$.  Furthermore, a graph $G$ of an arbitrary order $n\to\infty$ satisfies $\obj{G}=1/27+o(1)$
 if and only if $G$ has triangle density $o(1)$ and is almost $1/3$-regular. 
\end{theorem}

\bpf  Take any graph $G$.
For a vertex $u_i$ of $G$, we use the shortcuts $d_i:= \deg_G(u_i)$, $\Gamma_i:= \Gamma_G(u_i)$ and $\O\Gamma_i:= \Gamma_{\O G}(u_i)$.
Thus $\Gamma_i$ and $\O\Gamma_i$ are the sets of neighbours and non-neighbours respectively of $u_i$ in $G$ (which by definition do not contain the vertex $u_i$ itself). Also, for adjacent vertices $u_0$ and $u_1$, we  let $t = t(u_0, u_1)$ denote the number of triangles in $G$ containing the edge $u_0u_1$.
  
Let us show that $\Obj{G}\le n\,\pi_3(n-1)$ for any $n$-vertex graph $G$. Every embedding $f$ of $H$ into $G$ can be obtained by picking adjacent $u_0$ and $u_1$ for respectively $f(0)$ and $f(1)$, and then arbitrary $f(2)\in \Gamma_0\cap \O\Gamma_1$ and $f(3)\in \O\Gamma_0\cap \O\Gamma_1$; conversely, every such map $f$ gives an embedding.
 Since the sets from which we pick $f(2)$ and $f(3)$ partition the $(n-1-d_1)$-set $\O\Gamma_1$, the number of choices here is at most $\pi_2(n-1-d_1)$.
 Summing this over all ordered pairs $(u_0,u_1)$ of adjacent vertices, we obtain the desired upper bound
\begin{equation}
\label{eq:15ineq}
 \Obj{G}\le \sum_{u_1\in V(G)}\sum_{u_0\in \Gamma_1} \pi_2(n-1-d_1) = \sum_{u_1\in V(G)} d_1\pi_2 (n-1-d_1)\le n \, \pi_3(n-1),
\end{equation} 
 where we used the trivial observation that
\begin{equation*}
D( u ):= \pi_3(n-1)- \deg_G(u) \,\pi_2 (n-1- \deg_G(u) ),
\end{equation*}
 is non-negative for every vertex $u\in V(G)$. Let us call $D(u)$ the \emph{defect} of~$u$. 
 
For future reference,  we also define the \emph{defect} of adjacent vertices $u_0$ and $u_1$ to be
\begin{eqnarray*}
D(u_0,u_1)&:= &\pi_2(n-1-d_1)-| \Gamma_0\cap \O\Gamma_1| \cdot| \O\Gamma_0\cap \O\Gamma_1| \nonumber\\
   & \, =&\pi_2(n-1-d_1)-(d_0-1-t)(n-d_0-d_1+t)\ \ge\ 0,\label{eq:DefectFormula}
\end{eqnarray*}
 and the \emph{(total) defect} of~$G$ to be 
   \begin{equation*}
   D(G):= \sum_{u_1\in V(G)} \left(D(u_1)+\sum_{u_0\in \Gamma_1} D(u_0,u_1)\right)\ge 0.
\end{equation*}  
   Then the upper bound proved in~\eqref{eq:15ineq} is a consequence of the identity
\begin{equation}
\label{eq:15}
   \Obj[H_{15}]{G}=n\,\pi_3(n-1)-D(G).
\end{equation} 
  
Take any $n$-vertex graph $G$ attaining equality in~\eqref{eq:15ineq}; thus the defect of $G$ is zero, which implies that all vertex and pair defects are zero.  Take any adjacent $u_0$ and $u_1$. Since $D(u_0,u_1)=0$, it holds that
   \begin{equation*}
   -1\le (n-d_0-d_1+t)-(d_0-1-t)\le 1.
\end{equation*}
   By re-arranging the terms, we obtain
\begin{equation}
\label{eq:d0d1}
   n+2t\le 2d_0+d_1\le n+2t+2.
\end{equation}
   Also, the equality $D(u_1)=0$ (for any degree $d_1$) is equivalent to $| \,2d_1-(n-1-d_1)\,| \le 2$, which can be rewritten as
\begin{equation}
\label{eq:d1}
   n-3\le 3d_1\le n+1.
\end{equation}
   Note that \eqref{eq:d0d1} holds for every two adjacent vertices taken in either order, so it remains true if we swap $d_0$ and $d_1$. Thus it follows from~\eqref{eq:d0d1} and~\eqref{eq:d1} that $n+2t\le n+1$, that is, $t=0$. Since $u_0u_1\in E(G)$ was arbitrary, $G$ is triangle-free. 
   
  Let us show that $G$ is $\lceil n/3\rceil$-regular. Take any vertex $u_1$. By~\eqref{eq:d1}, $u_1$ cannot be an isolated vertex, hence we can pick a vertex $u_0$ adjacent to $u_1$. 
  If $d_0=d_1$ then the above equations imply that $n\le 3d_1\le n+2$ and thus $d_1=\lceil n/3\rceil$, as desired. If $d_0<d_1$ then, by \eqref{eq:d0d1} and~\eqref{eq:d1}, we have $n\le 2d_0+d_1\le 3d_1-2\le n-1$, a contradiction. If $d_1>d_0$ then we get a contradiction in the same way except we use 
  \eqref{eq:d0d1} with $d_0$ and $d_1$ swapped.
   
  Thus if $n$ is 0, 2, or 5 modulo 6 then every $n$-vertex graph attaining the upper bound is as stated (and conversely). Also, the set of such graphs is non-empty as shown by $B_{3k,2k}$ for $n=6k$, $B_{3k+1,2k+1}$ for $n=6k+2$, and $R_{3k-1,k}$ for $n=6k-1$ with integer $k\ge 1$.

  Next, let $n=6k+3$. Take any graph $G$ with $6k+3$ vertices and $\Obj{G}\ge n\,\pi_3(n-1)-2k$.
  Note that $D(u_1)=0$ if and only if $d_1$ is $2k$ or $2k+1$; otherwise it is easy to see that 
  \begin{equation*}
  D(u_1)\ge \pi_3(6k+2)-\max\left\{ (2k+2)(2k)^2,(2k-1)(2k+1)(2k+2)\right\}=2k
\end{equation*} 
  and equality holds if and only if $d_1=2k+2$. Thus every degree of $G$ is $2k$ or $2k+1$, except at most one vertex of degree~$2k+2$. 
  
  First, suppose that there is a vertex $u_1$ of degree $2k$. For any $u_0\in
  \Gamma_1$, we have that
\begin{equation*}
   D(u_0,u_1) + D(u_1, u_0)=\pi_2(4k+2)-(d_0-t-1)(4k+t+3-d_0) + \pi_2(6k+2-d_0)-(2k-t-1)(4k+t+3-d_0),
\end{equation*}
   which is 
   \begin{itemize}
   \item $2t^2 + 8t + 8$, at least 8 when $d_0=2k$,
   \item $2t^2 + 5t + 3$, at least 3 when $d_0=2k+1$,
   \item $2t^2 + 2t + 1$, at least 1 when $d_0=2k+2$,
   \end{itemize}
with equality in these cases only if $t=0$. Since at most one vertex has degree $2k+2$, we have that 
$$\sum_{u_0\in \Gamma_1} \left(D(u_0,u_1) + D(u_1, u_0)\right) \ge 1 + 3(2k-1)= 6k-2.
$$
Hence, when we have a vertex with degree $2k$, it is impossible to have $D(G)\le 2k$.
   
  Second, suppose that each degree is $2k+1$ or $2k+2$. Since $n$ is odd, there is at least one vertex $u_1$ of degree $2k+2$. We have
   \begin{equation*}
   D(2k+2)\ge \pi_3(6k+2)-(2k+2)\pi_2(4k)=2k\cdot (2k+1)^2 - (2k+2) (2k)^2=2k.
\end{equation*}
   Thus all other defects are zero, $u_1$ is the only vertex of degree $2k+2$, all other vertices have degree $2k+1$ and there are no triangles, as desired.

  Conversely, the above calculations show that, for each $k\ge 1$, every graph $G$ as stated in the theorem satisfies $\Obj{G}=n\,\pi_3(n-1)-2k$.
   
For a construction of at least one such graph, we start with $R_{3k+1,k}$. We will add edges of the form $\{i,i+2k\}$ or $\{i,i+2k+1\}$ for some residue $i$. No triangle can contain exactly one of these added edges: each original edge of $R_{3k+1,k}$ has cyclic length at least $2k+2$ so after two jumps in the graph $R_{3k+1,k}$ (in the same or opposite directions) we are at distance at most $2k-1$ from where we started but no edges of such length are added. Note that the graph on residues modulo $6k+3$ consisting only of the edges of length $2k$ is a vertex-disjoint union of cycles. Call this graph $H$. The graph $H \setminus \{0, 2k+1, 4k+2\}$ has a perfect matching. Indeed, if $3$ divides $k$, then $H$ is the union of three cycles, each of length $2k+1$, and the three removed vertices lie in different cycles; deleting them leaves three paths on $2k$ vertices. If $3$ does not divide $k$, then $H$ is a single cycle of length $6k+3$, and the vertices $\{0,2k+1,4k+2\}$ are equally spaced on this cycle; deleting them leaves three paths, again each on $2k$ vertices. In either case, the result has a perfect matching $M$ of size $3k$.

Let $G$ be obtained from $R_{3k+1,k}$ by adding the two edges
$\{0,2k+1\}, \{2k+1,4k+2\}$ and the edges from $M$. The only two edges sharing a vertex are $\{0,2k+1\}, \{2k+1,4k+2\}$. The third edge needed to complete a triangle on these three vertices would be $\{0,4k+2\}$, which has cyclic length $2k+1$. The other added edges are vertex-disjoint, hence this graph is triangle-free and has the required degree sequence.
  
Let $n=6k+1$ with integer $k\to\infty$.
Take any order-$n$ graph $G$. First, let us show that the defect of $G$ is at least $6k^2+O(k)$. Suppose that, say, $D(G)\le 6k^2$ as otherwise we are done.
  
  Note that $D(u_1)=(2k)^3-d_1\pi_2(6k-d_1)$. It is 0 if $d_1=2k$, it is $2k$ if $d_1=2k\pm1$ and, for $d_1=2k+s$ with $| s| \ge 2$, it is at least
\begin{equation}
\label{eq:Du1}
   (2k)^3-(2k+s)\,\left(\frac{4k-s}2\right)^2=\frac{3}2\, ks^2 - \frac14\, s^3\ge \max\left\{6k-2,\frac{1}2\, ks^2\right\}
\end{equation}
   since $| s| \le 4k$. Let $A$ (resp.\ $B$) be the set of vertices of degree $2k$ (resp.\ $2k+1$). Let $a:= | A| $ and $b:= | B| $. By above, we have that $(n-a)\cdot 2k\le D(G)\le 6k^2$. Thus $a\ge n/2$.
   
  Note that for any two adjacent vertices $u_0,u_1\in A$, it holds that
\begin{equation}
\label{eq:6k+1D01}
   D(u_0,u_1)=\pi_2(4k)-(2k-1-t)(2k+1+t)=(t+1)^2 \ge 1.
\end{equation} 
   Thus every edge inside $A$ contributes at least 2 to $D(G)$, since it is counted for both its orderings. By $D(G)\le 6k^2$ and~\eqref{eq:Du1}, we can assume that the maximum degree is at most $2k+\sqrt{12 k}$. Thus the number of edges between $A$ and its complement is at most $b(2k+1)+(n-a-b)(2k+\sqrt{12 k})$ and 2 times the number of edges inside $A$ can be lower bounded by $a\cdot 2k$ minus this quantity.
   By~\eqref{eq:Du1} and~\eqref{eq:6k+1D01}, we have
\begin{equation}
\label{eq:6k+1D}
   D(G)\ge \left(a\cdot 2k-b(2k+1)-(n-a-b)(2k+\sqrt{12 k})\right)+b\cdot 2k + (n-a-b)\cdot (6k-2).
\end{equation}
   This is a linear function of $b$ with the coefficient at $b$ being $-4k+\sqrt{12k}+1<0$. Thus it is at least its value when $b=n-a$, giving $D(G)\ge 2ak+O(k)\ge 6k^2+O(k)$, as desired.
   
The stated lower bound on $\Obj{6k+1}$ comes from the following construction. Let $a:= \lfloor k/2\rfloor$. We start with $B_{3k,2k}$, denoting its parts as $U:= \{u_0,\dots,u_{3k-1}\}$ and $W:= \{w_0,\dots,w_{3k-1}\}$ where index arithmetic will be taken  modulo $3k$. It will be notationally convenient to rotate one set of the residues, so we connect $v_i$ to $u_j$ in $B_{3k,k}$ if and only if $i+j$ belongs to the interval $\{2k+1,\dots,4k\}$ modulo $3k$. Add new $a$ edges $\{u_i,w_{i+k-a}\}$ for $i\in \{0,\ldots,a-1\}$. Add a new vertex $u$ connected to $2k$ elements $w_{k},\dots,w_{3k-1}$. The obtained graph $G$ is clearly bipartite. It has $2k+2a$ vertices of degree $2k+1$ (namely, $u_0,\dots,u_{a-1}$ and $w_{k-a},\dots,w_{3k-1}$) with the remaining $4k+1-2a$ vertices (namely, $u_{a},\dots,u_{3k-1}$, $w_0,\dots,w_{k-1-a}$ and $u$) having degree $2k$, so $\sum_{v\in V(G)} D(v)=(2k+2a)\cdot 2k$. It is routine to see that, since there are no triangles, if we take two adjacent vertices, one of degree $2k+1$ and the other of degree $2k+1$ or $2k$ then their defect (in either ordering) is zero. Also, as we observed earlier (see~\eqref{eq:6k+1D01}) the defect of any ordered adjacent pair of vertices each of degree $2k$  is~$1$. Each such pair contains one of $w_0,\dots,w_{k-a-1}$, while each of these vertices $w_i$ is adjacent to all degree-$(2k+1)$ vertices in~$U$ and thus has  exactly $2k-a$ neighbours of degree $2k$. Also, all $2k$ neighbours of the added vertex $u$ have degree $2k+1$. 
  Thus the total defect  is
  $(2k+2a)\cdot 2k + 2\cdot (k-a)(2k-a)=15k^2/2+O(k)$, as desired.

  Finally, let $n=6k+4$. Here the proof of the upper bound is very similar to the case $n=6k+1$, so we will be rather brief. Take any extremal graph $G$. For $u_1\in V(G)$, it holds that $D(u_1)=0$ only if $d_1=2k+1$. Furthermore, $D(u_1)=2k+1$ if $d_1=2k$ or $2k+2$ and then $D(u_1)$ grows quadratically as $d_1$ moves away from $2k+1$. In particular, the maximum degree of $G$ is at most $2k+O(\sqrt{k})$. Also, one can show that every adjacent pair of vertices, each of degree $2k+1$, has defect at least 1. It follows from the above observations by an easy adaptation of the argument leading to~\eqref{eq:6k+1D} that  $D(G)\ge 6k^2+O(k)$.
  
  For the lower bound, observe that, if the graph is triangle-free, then the defect of any adjacent pair of vertices (in either order), one of degree $2k+2$ and the other of degree $2k+1$ or $2k+2$, is zero. Our construction will have only vertices of degree $2k+1$ and $2k+2$, about $n/2$ of each degree; then we have to minimise the number of edges spanned by a set of vertices of degree $2k+1$. Let $a:= \lceil (3k+2)/2\rceil$. Take the union of  complete bipartite graphs $\MultiPartite{a,\,3k+2-a}$ with parts $U_0\cup U_1$ and $\MultiPartite{3k+2-a,\,a}$ with parts $W_0\cup W_1$. Add a copy of $B_{a,a-k}$ with parts $U_0$ and $W_1$ and a copy of $B_{3k+2-a,2k+1-a}$ with parts $U_1$ and $W_0$. All vertices in $U_0\cup W_1$ (resp.\ $U_1\cup W_0$) have degree $2k+2$ (resp.\ $2k+1$). Also, each vertex of degree $2k+1$ has $k/2+O(1)$ neighbours of degree $2k+1$. Thus the total defect is
   \begin{equation*}
   | U_0\cup W_1| \cdot k/2 + | U_1\cup W_0| \cdot (2k+1) +O(k) = \frac{15}2\, k^2+O(k),
\end{equation*}
  as desired.

Finally, it is routine to see from~\eqref{eq:15} that a graph $G$ of an arbitrary order $n\to\infty$ is almost extremal
 (that is, $D(G)=o(n^4)$) if and only if it is as stated in the theorem.  
\epf
  
The upper bound of Theorem \ref{th:H15Exact} for $n=6k+1$ and $n=6k+4$ can be improved  with extra work to $\Obj{n}\le n\,\pi_3(n-1)- (6+c+o(1))k^2$ for some (small) constant $c>0$. However, it seems that pushing $6+c$ all the way up to $15/2$ (if this is the correct constant at $k^2$) would be rather messy. So we limit ourselves to the stated bounds which are relatively easy to show.

\subsection{Computer-generated results for  $H_9$, $H_{10}$ and $H_{11}$}\label{se:computer}

Here we present the cases (namely the graphs $H_9$, $H_{10}$ and $H_{11}$) where the flag algebras not only give the semi-inducibility constant but also prove perfect stability by applying the criterion given by the second author, Slia\v can and Tyros~\PST{Theorem 7.1} that can be automatically checked by computer. 

Since the statement of~\PST{Theorem 7.1} is rather long, we just give an informal description of it. In addition to a flag algebra certificate proving $\obj[\gamma]{}\le u$ as in~\eqref{eq:FAMain} (with the matching lower bound coming from blowups of some base graph $B$), we have to provide some graph $\tau$ with at most $\N-2$ vertices that satisfies certain properties. Some of the required properties of $\tau$ are:
(a) if we forbid $\tau$ as induced subgraph then the maximum value of $\obj[\gamma]{}$ on graphs of order $n\to\infty$ drops by $\Omega(1)$, (b) $\tau$ admits the unique (up to an automorphism of $B$) homomorphism $f:\tau\to B$, and (c) each vertex of $B$ is uniquely identified by its adjacencies to $f(\tau)$. The gist of the proof of ~\PST{Theorem 7.1} is that every almost extremal $G$ of order $n\to\infty$ has $\tau$-density $\Omega(1)$ by (a), every copy of $\tau$ in $G$ gives pairwise disjoint $V_0,\dots,V_{m-1}\subseteq V(G)$ naturally defined by looking at adjacencies of a vertex of $G$ to that copy  by (b)--(c) so that any single discrepancy between $E(G)$ and $E(\blow{B}{V_0,\dots,V_{m-1}})$ gives $\Omega(n^{N-v(\tau)-2})$ $\N$-vertex subgraphs $F$ of $G$ that are `wrong' (i.e.\ do not admit a homomorphism to $B$). If, furthermore, every `wrong' $F$ has positive slack $c_F>0$ then $G$ must be close to $\blow{B}{V_0,\dots,V_{m-1}}$  for a `typical' copy of $\tau$. The proof also needs that at least one of certain Conditions (i) and (ii) of~\PST{Theorem 7.1} holds, and if Condition~(i) holds then the proof gives that the optimal $B$-vector is unique up to the automorphisms of $B$. In brief, Condition~(i) states that the $\mathrm{SOS}$-term in~\eqref{eq:FAMain} contains squares of $\tau$-rooted subgraph densities  such that the positive semi-definite matrix that encodes them has the smallest possible co-rank 1. When we use~\PST{Theorem 7.1}, we just list $\N$ and $\tau$, and state which of Conditions~(i) and (ii) was verified by our code.

\newcommand{\SIRegularTheoremTemplate}[8]{Let $H_{#8}$ be the graph on $[4]$ with the red edge set $\{#1\}$ and the blue set $\{#2\}$. Then $\obj[H_{#8}]{}=#3$, where the upper bound is proved via flag algebras with $\N=#5$ while the lower bound comes from #7. Also, Corollary~\ref{co:reg} applies here and gives that every almost extremal sequence of graphs is almost $a$-regular (and an easy calculation shows that necessarily $a=#6$).}

\ExplSyntaxOn
\NewDocumentCommand{\CleanString}{m}
 {
  \tl_set:Nn \l_tmpa_tl { #1 }
  \tl_replace_all:Nnn \l_tmpa_tl { / } { }
  \tl_replace_all:Nnn \l_tmpa_tl { , } { ,~ }
  \tl_use:N \l_tmpa_tl
 }
\ExplSyntaxOff

\newcommand{\SemindStableTheoremTemplate}[9]{Let \[H := H_{#9} = \DrawSemiGraphLabeled{4}{#1}{#2}{0,1,2,3},\] be the graph on $[4]$ with the red edge set $\{\CleanString{#1}\}$ and the blue edge set $\{\CleanString{#2}\}$. Then it holds that $\obj[H_{#9}]{}=#3$, the problem is perfectly $B$-stable for $B = #4$ by \PST{Theorem~7.1(#8)} with $\N=#5$ and $\tau=#6$. Furthermore, the unique maximizer of $\obj[H_{#9}]{}\left(\blow{B}{\V a}\right)$ is $\V a = \left( #7 \right)$.
}

\begin{theorem}\label{th:semiind9}
\SemindStableTheoremTemplate{0/1,0/2,1/2,2/3}{0/3}{12/125}{K_5}{7}{([5], \{01\})}{1/5, 1/5, 1/5, 1/5, 1/5}{i}{9}\qed 
\end{theorem}

\begin{theorem}\label{th:semiind10}
\SemindStableTheoremTemplate{0/1,1/2}{0/2,0/3}{1/8}{K_2}{5}{K_1}{1/2, 1/2}{ii}{10}
\end{theorem}

\bpf Since Condition~(i) of \PST{Theorem~7.1} does not apply here, we have to check the uniqueness of the maximiser~$\V a$.
Note that $\Obj[H]{\MultiPartite{xn,\,(1-x)n}}=p(x) n^{\fallingfact{4}}+o(n^4)$, where $p(x):=x(1-x)^3+x^3(1-x)$ since the only part assignment that gives a copy of $H$ is when the vertices $0,2,3$ are in one part and the vertex $1$ in the other. As it is easy to see,
the maximum of $p(x)$ on $[0,1]$ is $1/8$ and this is attained if and only if $x=1/2$, as desired.  
\epf

\begin{theorem}\label{th:semiind11}
\SemindStableTheoremTemplate{0/2,0/3}{0/1,2/3}{1/8}{K_2}{4}{K_2}{1/2, 1/2}{i}{11}\qed
\end{theorem}

Since the flag algebra proof of $\obj[H_{11}]{}\le 1/8$ uses only $\N=4$, the identity in~\eqref{eq:FAMain} can be explicitly written. Namely, after some cleaning up when we try to make as many coefficients inside the certificate as possible to be 0,  the proof amounts to checking that, for any order-$n$ graph $G$, if we expand
 \begin{equation}\label{eq:H11SOS}
 \frac{n^4}8-\Obj[H_{11}]{G}-
 \frac38 \sum_{uw\in E(G)} (\deg_G(u)-\deg_G(w))^2 
 -\frac18 \sum_{uw\in E(\O G)} \left(|\Gamma_G(u)\cap \Gamma_G(w)|-|\Gamma_{\O G}(u)\cap\Gamma_{\O G}(w)|\right)^2
 \end{equation}
as a linear combination of 4-vertex subgraph counts up to error term $o(n^4)$ then each coefficient will be non-negative. While this proves that $\obj[H_{11}]{}\le 1/8$, it does not give much insight. (Unfortunately, computer calculations indicate that there is no proof of $\obj[H_{11}]{}\le 1/8$ with $\N=4$ as in~\eqref{eq:H11SOS} that  has only one sum, either over edges or over non-edges of $G$.)

\subsection{Graph $H_6$}\label{se:H6}

For the graph $H_6$, flag algebras give the sharp upper bound on the semi-inducibility constant. However, perfect stability does not hold. So we prove the Erd\H os-Stone stability by extracting some further information from the certificate, in addition to just the upper bound $\obj[H_6]{}\le 1/8$.
Then, by using the standard stability approach, we get exact result for all $n\ge n_0$.

\begin{theorem}\label{th:semiind6}
Let \[H_6 = \DrawSemiGraphLabeled{4}{1/2}{0/2,0/3,2/3}{1,2,0,3}\] be the graph on $[4]$ with the red edge set $\{02\}$ and the blue set $\{01,03,13\}$. Then $\obj[H_{6}]{}=1/8$, where the upper bound is proved via flag algebras with $\N=5$. The lower bound comes from blowups of $K_2$ with the part ratios given by $\V a=(1/2, 1/2)$.
Moreover, if a 5-vertex graph $F$ has density $\Omega(1)$ in an almost extremal graph then $F$ is complete bipartite or is obtained from a 4-vertex complete bipartite graph by adding one isolated vertex.
\end{theorem}
\bpf The result is obtained by a standard application of flag algebras. Also, our script checks that every  $5$-vertex graph $F$ with  $c_F=0$ is complete bipartite plus at most one isolated vertex.\epf

First, let us observe that the semi-inducibility problem for $H_6$ does not satisfy perfect stability. Consider the graph $G:=\MultiPartite{x_0 n,\, x_1 n}+\O{K}_{(1-x_0-x_1)n}$ (that is, $G$ is obtained from the complete bipartite graph with part sizes $x_0 n$ and $x_1 n$ by adding $(1-x_0-x_1)n$ isolated vertices). We have that $\obj[H_6]{G}=q(x_0,x_1)+o(1)$, where 
\begin{equation}\label{eq:H6Constr}
q(x_0,x_1):= x_0x_1\cdot  \left( 1-2x_0x_1+(1-x_0-x_1)^2\right).
\end{equation}
 Indeed, we first choose an unordered pair $uw$ of adjacent vertices in $G$ (to be the image $\{f(0),f(2)\}$); then an ordered pair $(f(1),f(3))$ of two non-adjacent vertices of $G$ gives rise to exactly two possible maps of $02$ to $uw$ if both $f(1)$ and $f(3)$ are disconnected from $f(0)$ and $f(2)$, otherwise exactly one map. 
 Now, if we take small $\e>0$ and let each of $x_0$ and $x_1$ be $1/2-\e$ then the density of $H_6$ is $1/8-\e^2+O(\e^3)$ while $G$ can be shown to be at the edit distance at least $\Omega(\e n^2)$ from a blowup of~$K_2$. Thus there is no perfect stability.

On the other hand, we can show that this problem satisfies the Erd\H os-Simonovits stability.

\begin{lemma}\label{lm:H6Stab} 
 For every $\e>0$ there are $\delta>0$ and $n_0$ such that every graph $G$ with $n\ge n_0$ that satisfies $\obj[H_6]{G}\ge 1/8-\delta$ is within edit distance at most $\e n^2$ from $\MultiPartite{n/2,\,n/2}$.\end{lemma}
 \bpf
 Suppose that the lemma is false for some $\e>0$. Thus we can find  an almost extremal graph of order $n\to\infty$ that is $\e n^2$ far from $\MultiPartite{n/2,\,n/2}$. By the second part of Theorem~\ref{th:semiind6} and the Induced Removal Lemma~\cite{AlonFischerKrivelevichSzegedy00}, we can change $o(n^2)$ adjacencies in $G$ such that every 5-vertex subset spans a complete bipartite graph apart from possibly one isolated vertex. In particular, the new graph $G$ contains neither a triangle nor an induced path on 4 vertices. It follows that $G$ is bipartite, for a shortest odd cycle in $G$, which is necessarily induced, has to contain one of these two configurations. Furthermore, if two vertices $u,w$ are connected by a path of length 3 then $u,w$ must be adjacent. It follows that each component of $G$ is a complete bipartite graph. Since $G$ has no induced matching with 2 edges, it is a complete bipartite graph $\MultiPartite{n_0,\,n_1}$ plus $n-n_0-n_1$ isolated vertices. Let $x_0:=n_0/n$ and $x_1:=n_1/n$. Thus $\obj[H_6]{G}=q(x_0,x_1)$, where $q$ is the degree-4 polynomial defined by~\eqref{eq:H6Constr}. 

It follows from Theorem~\ref{th:semiind6} that $q(x,y)\le 1/8$ on 
\[
S:=\{(x,y)\in \I R^2\mid x\ge 0,\ y\ge 0,\ x+y\le 1\}.
\] 
We need to show that equality holds only if $x=y=1/2$. For this task, we first re-prove the upper bound $1/8$ on the values of $q$ on $S$ using basic calculus. Here we maximise a polynomial on a compact set, so the maximum is attained at a critical point or on the boundary. First, let us find all critical points $(x,y)$ of $q$ in the interior of $S$. If $y=x$ then
\[
 q(x,x)=\frac18+\frac18(2x-1)^2(4x^2-4x-1)
\]
 is less than $1/8$ for each $x\in (0,1/2)$. So assume that $x\not=y$. The following identity is routine to check:
 $$
 \frac{\partial q(x,y)}{\partial y}-\frac{\partial q(x,y)}{\partial x}=(x-y)\left((x-y)^2+2(1-x-y)\right).
 $$
 Its right-hand side is non-zero so $(x,y)$ is not a critical point. Thus it remains to consider the boundary of~$S$. For every real $x$, 
we have $q(x,0)=q(0,x)=0$ while $q(x,1-x)=1/8-(1-2x)^2/8$ is at most $1/8$ with equality if and only if $x=1/2$. Thus $(1/2,1/2)$ is the only point of $S$ where $q$ attains value~$1/8$.

By the continuity of the involved functions, it follows that each of $x_0$ and $x_1$ is $1/2+o(1)$. However, then $G$ is within edit distance $o(n^2)$ from $\MultiPartite{n/2,\,n/2}$, contradicting our assumption.\epf

Although perfect stability does not hold for $H_6$, we will show  in Theorem~\ref{th:H6Exact} that every extremal graph $G$ of sufficiently large order  is exactly complete bipartite.  Our proof will need the following lemma which is routine to verify; in fact, one can explicitly write all optimal $a$ if needed, see e.g.~\cite[Eq.~(16)]{BodnarPikhurko25ind}.

\begin{lemma}\label{lm:Opt3+1} As $n\to\infty$, an integer $a$ with $0\le a\le n$ that maximises the polynomial 
 \begin{equation}\label{eq:P}
 P_n(x):=x^{\fallingfact{3}}(n-x)+x(n-x)^{\fallingfact{3}}.
 \end{equation}
 satisfies $a=n/2+\sqrt{3n}/2+O(1)$ and the maximum value is $n^{\fallingfact{4}}/8 +3n^2/4+O(n)$.\qed
 \end{lemma}

\begin{theorem}\label{th:H6Exact} Let $H$ be the graph $H_6$ (that has one red edge $02$ and blue edges $01,03,13$). Then every graph $G$ of sufficiently large order $n$ with $\Obj[H]{G}=\Obj[H]{n}$ is a blowup of $K_2$.\end{theorem}
 
\bpf Let $G=(V,E)$ be any extremal graph of sufficiently large order $n$. For notational convenience, we view $n$ as going to infinity and use the asymptotic notation accordingly.
Let $V_0\cup V_1=V$ be a max-cut partition of $G$. By Lemma~\ref{lm:H6Stab}, we have that the number of edges between $V_0$ and $V_1$ is at least $(1/4+o(1))n^2$. It follows that $|V_i|=n/2+o(n)$ for each $i=0,1$ and almost every pair between $V_0$ and $V_1$ is an edge of~$G$. Also, it follows that the number of edges inside $V_0$ or $V_1$ is $o(n^2)$ as otherwise $G$ will be $\Omega(n^2)$-far from being bipartite (e.g.\ by having $\Omega(n^3)$ triangles). Let $W$ consist of all \emph{wrong} pairs, that is, non-edges across the parts and edges inside a part. As we just argued, it holds that $|W|=o(n^2)$.

Take any vertex $u\in V$. For $i=0,1$, define $y_i:=|\Gamma(u)\cap V_i|/|V_i|$ to be the fraction of vertices in $V_i$ that are neighbours of~$u$. 

Let us show that $\Obj[H]{G,u}$, the number of embeddings of $H$ into $G$ that use the vertex $u$, is equal to $p(y_0,y_1)n^3+o(n^3)$, where we define
\beq{eq:p}
 p(x,y):=\frac{x^3+x^2 y+x y^2+y^3}8- \frac{x y}2-\frac{x}{8}-\frac{y}{8}+\frac{1}{2}.
 \eeq
  This calculation is straightforward but rather lengthy, so we include some details to make it easier to verify for the reader. Here, we can replace $G-u$ by the blowup of $K_2$ with the parts $V_0\setminus\{u\}$ and $V_1\setminus\{u\}$, since this changes the number of embeddings via $u$ by at most $4!\, |W|\, n=o(n^3)$.
  Clearly, $\Obj[H]{G,u}$ is the sum over all $i$ in $V(H)=[4]$ of the number of embeddings that map $i$ to $u$. Also, the $i$-th summand can be written as $p_i(y_0,y_1)n^3+o(n^3)$, where
  \begin{eqnarray*}
  p_{0}(x,y)&:=& \frac{x+y}2\cdot \frac{(1-x)^2+(1-y)^2}{4},\\
  p_{1}(x,y),\, p_{3}(x,y)&:=& \frac{(1-x)^2}8+\frac{(1-y)^2}8,\\
  p_{2}(x,y)&:=&\frac{x+y}8.
  \end{eqnarray*}
  For example, the formula for $p_{1}$ follows by observing that 
$0$ and $3$ have to be mapped inside the same part $V_i$ and be non-neighbours of $u$, while $2$ has to go to a vertex in the other part $V_{1-i}$ (and, conversely, every such injection into $V\setminus\{u\}$ gives an embedding of $H$). Finally, it is routine to check that $p=\sum_{i=0}^3 p_i$. 

The pairs $(x,y)=(0,1)$ and $(1,0)$ correspond to $u$ being adjacent to every vertex in one part and to none in the other. As expected, the value of $p$ on them is $1/2=4\,\obj[H]{}$. Also, it holds that $p(x,y)\le 1/2$ for every $(x,y)\in [0,1]^2$ (as otherwise we would contradict Theorem~\ref{th:semiind6} by  $\e n$ vertices to $G$ with appropriate attachment for some small constant $\e>0$). Recall that the strictness property would state here that all optimal $(x,y)$ are of this form. However, this property does not hold for  $H$-semi-inducibility, since $p$ attains the maximum value $1/2$ also on $(0,0)$. As a step in our exactness result, we need to describe all optimal $(x,y)$. 
For this, we have to calculate the maximum of the explicit polynomial $p$ from scratch, using basic calculus. 

\begin{claim}\label{cl:p}
 The maximum of the polynomial $p(x,y)$ (defined in~\eqref{eq:p}) over $(x,y)\in [0,1]^2$ is $1/2$ and this is attained if and only if $(x,y)$ is one of $(0,1)$, $(1,0)$ and $(0,0)$.\end{claim}
\bpf[Proof of Claim~\ref{cl:p}.] We maximise a differentiable function over a compact set so its maximum is attained at a critical point in the interior or at some point on the boundary. The identity 
 \[
 \frac{\partial p(x,y)}{\partial x}-\frac{\partial p(x,y)}{\partial y}=\frac{(x-y)(x+y+2)
}4 \]
 shows that all interior critical points are limited to the line $x=y$. We have that $p(x,x)-1/2=x(2x^2-2x-1)/4$ and this is negative for all $0<x<1$. So it remains to consider the boundary. For $0\le x\le 1$, it holds that $p(x,0)-1/2= x^3/8-x/8\le 0$  with equality if and only if $x=0$ or $x=1$, and that $p(x,1)-1/2=x(x^2+x-4)/8\le 0$ with equality if and only if $x=0$. This takes care also of the remaining pieces of the boundary by the symmetry between $x$ and $y$, finishing the proof.\epf

The average number of embeddings that use a random vertex of the extremal graph $G$ is $\Obj[H]{n}\cdot 4/n$. Since we cannot increase $\Obj[H]{G}$ by replacing a vertex by a clone of another,  every vertex $u$ is in at least $4\,\obj[H]{n}n^3+O(n^2)=(4\,\obj[H]{}+o(1))n^3$ embeddings. This implies by Claim~\ref{cl:p} that the vector $(y_0,y_1)$ (where $y_i:=|\Gamma(u)\cap V_i|/|V_i|$) is $o(1)$-close to $(0,1)$, $(1,0)$, or $(0,0)$. In the first two cases, the vertex $u$ belongs necessarily to respectively $V_0$ and~$V_1$ by the max-cut property.

Define $V_2$ to be the set of those vertices $u\in V$ for which the last alternative holds (which is equivalent to $\deg(u)=o(n)$). The edit distance from $G$ to a complete partite graph is $\Omega(|V_2|\, n + |W|)$. Hence, by stability (that is, by Lemma~\ref{lm:H6Stab}), we have that $|V_2|=o(n)$. 

Let $V'_i:=V_i\setminus V_2$ for $i=0,1$. 
Define $W'$ to consist of the non-edges between $V'_0$ and $V'_1$, the edges inside $V'_0$ or $V'_1$ (so far these the same definitions same as for $W$) and the edges between $V'_0\cup V'_1$ and~$V_2$. By above, we have that the maximum degree $\Delta(W')=o(n)$. As the next step, we can show that $W'$ is in fact empty.
 
\begin{claim}\label{cl:W'} It holds that $W'=\emptyset$.\end{claim}

\bpf[Proof of Claim~\ref{cl:W'}.] Suppose on the contrary that some pair $uw$ is in~$W'$. Let $G':=G\oplus uw$ be obtained from $G$ by changing the adjacency between $u$ and $w$.  We would like to show that $\Obj[H]{G'}-\Obj[H]{G}\ge \e n^2$ for some small constant $\e>0$, which will be the desired contradiction. Note that only those injections that use both $u$ and $w$ can contribute non-zero amount to  $\Obj[H]{G}-\Obj[H]{G'}$.
Thus if, in each of $G$ and $G'$, we remove all vertices of $V_2\setminus \{u,w\}$ and correct every pair of $W'$ (or, equivalently, $W$) except $uw$, then the difference $\Obj[H]{G'}-\Obj[H]{G}$ changes by at most $4!\,(|V_2|n+|W'|+\Delta(W')n)=o(n^2)$, which can be taken care of by just slightly decreasing $\e$. Thus we can calculate $\Obj[H]{G'}-\Obj[H]{G}$ with respect to the simplified graphs.

First, suppose that $uw$ is an edge inside, say, $V'_0$. In this case, $G'$ is obtained from $G$ by making $u$ and $w$ non-adjacent.
Each map $f$ that is an embedding of $H$ into $G$ but not into $G'$ has to map $02$ to $uw$ in some order (2 choices). Ignoring $V_2$, every extension has to map  both $1$ and $3$ to $V'_0$, so we have in total at most $2 (n/2)^2+o(n^2)$ such embeddings. On the other hand, we can lower bound the number of the embeddings created by flipping $uw$ as follows. First, we can map $\{1,3\}$ to $uw$ (2 choices) then $0$ 
to $V'_0\setminus\{u,w\}$ and $2$ to 
$V'_1$ ($(n/2)^2+o(n^2)$ choices). Second, we can map $0$ into $uw$ and then one of $1,3$ to the remaining vertex in $uw$ (4 choices), then we can map $2$ to $V'_1$ and the remaining vertex of $\{1,3\}$ to $V'_0\setminus\{u,w\}$, having here $(n/2)^2+o(n^2)$ choices. Thus 
 \[
 0\ge \Obj{G'}-\Obj{G}\ge 2\cdot n^2/4+4\cdot n^2/4-2\cdot n^2/2+o(n^2)=n^2/2+o(n^2)>0,
 \] 
 a contradiction showing that $V'_0$ (and similarly $V'_1$) spans no edges. 
 
 Second, suppose that $(u,w)\in V'_0\times V'_1$ is a non-edge. Here $G'$ is obtained from $G$ by changing $uw$ to an edge. There are no maps that embed $H$ into $G$ but not into $G'$ as each has to map a blue edge to $uw$ leaving no choice for the third vertex of the blue triangle of $H$.  On the other hand, the number of maps which are embeddings into $G'$ but not $G$ is at least $2\cdot (n/2)^2+o(n^2)$: map $\{0,2\}$ to $uw$ in some order and then map $1$ and $3$ into the same part as $0$. This contradiction shows that all pairs between $V'_0$ and $V'_1$ are edges.
 
 Third, suppose that $u\in V_2$ and, say, $w\in V'_0$ (and these are adjacent in~$G$).
 Maps that are embeddings into $G$ but not $G'$ have to map $\{0,2\}$ to $uw$. If $0$ is mapped to $w$ then $1,3$ have both to be mapped to $V'_0$ giving at most $(n/2)^2+o(n^2)$ choices; if $0$ is mapped to $u$ then $1$ and $3$ have to be mapped into a same part (any) so we have at most $n^2/2+o(n^2)$ choices here.
 On the other hand, all the following injections give embeddings into $G'$ but not into $G$. We can map $1,3$ to $uw$ (2 choices), send $0$ to $V'_0\setminus\{w\}$ and send $2$ to $V'_1$ ($(n/2)^2+o(n^2)$ choices). Or we can map $3$ to $u$, $0$ to $w$, $1$ to $V'_0\setminus\{w\}$ and $2$ to $V'_1$ ($n^2/4+o(n^2)$ choices); we obtain the same bound when we map $1$ to $u$ and $0$ to $w$. Thus
  $$
  0\ge \Obj{G'}-\Obj{G}\ge 2(n^2/4) + 2 (n^2/4)-(3/4)n^2+o(n^2)=n^2/4+o(n^2)>0,
  $$
  a contradiction. This proves the claim.\epf

 Denote $n_i:=|V'_i|$ for $i=0,1$ and $n_2:=|V_2|$.
 As our next intermediate step, let us show that $n_2\le 30$. We do not optimise the bound  since we will later show that $n_2=0$ (by using that $n_2=O(1)$).
 
 Observe that by Lemma~\ref{lm:Opt3+1}
 we have
  \[
  \Obj{n}\ge \max_{a\in [n+1]} \Obj{\MultiPartite{a,\,n-a}}
  =\frac 18 n^{\fallingfact{4}} +\frac34\, n^2+O(n),
  \]
  and the average number of embeddings per vertex in extremal $G$ is at least $(4/n)\Obj{n}=(n-1)^{\fallingfact{3}}/2+O(n)$. Since we cannot increase $G$ by removing a vertex and adding a clone of another vertex, it holds  for every vertex $u$  that the number of embeddings that use $u$ is at least 
  \begin{equation}\label{eq:Gu}
   \Obj[H]{G,u}\ge \frac{(n-1)^{\fallingfact{3}}}2-4!\, {n\choose 2}+O(n)=\frac{n^3}2-15n^2+O(n).
  \end{equation}
  Let us estimate $\Obj[H]{G,u}$ for any given $u\in V_2$. First, the number of the embeddings $f$ where the 3 other vertices of $H$ are all mapped into $V'_0\cup V'_1$ is at most $2(n_0^{\fallingfact{2}}n_1+n_0n_1^{\fallingfact{2}})$. Indeed, every such $f$ has to map $1$ or $3$ to $u$ and if, say, $f(1)=u$ then  $f(0)$ and $f(3)$ have to be in the same part $V'_i$ while $2$ is in the other part $V'_{1-i}$ for some $i\in \{0,1\}$. Next, let us show that the number of embeddings that use exactly one other vertex of $V_2$ is at most $|V_2|(n^2+o(n^2))$. We have to sum over all choices of the other vertex $w\in V_2\setminus \{u\}$.  If $uw$ is an edge then $0,2$ have to be mapped to $u,w$, and then $1,3$ are mapped to the same part $V'_i$ (for $i=0,1$); thus we have at most $2n^2/2+o(n^2)=n^2+o(n^2)$ embeddings for such $w$. If $uw$ is not an edge, then $1,3$  have to be mapped to $u,w$, and then $0,2$ are mapped to different parts among $V'_0$ and $V'_1$; again we get at most $n^2+o(n^2)$ embeddings and the claim follows. 
  Finally, the number of embeddings that use at least two other vertices from $V_2$ is at most $O(n_2^2n)=o(n_2n^2)$. Thus by~\eqref{eq:Gu} we have
  \begin{eqnarray*}
  0
    &\ge& \left(\frac{n^{3}}2-15n^2\right)-\left(2\left(n_0^{\fallingfact{2}}n_1+n_0n_1^{\fallingfact{2}}\right)+n_2n^2\right)+o(n_2n^2)\\
   &=& \frac{(n_0-n_1)^2(n_0+n_1)}2+4n_0n_1 + \frac{(n_0+n_1)^2(n_2-30)}2+o(n_2n^2),
  \end{eqnarray*}
  where the last identity is obtained by substituting $n=n_0+n_1+n_2$. This cannot be satisfied for any $n_2\ge 31$, proving the claimed upper bound on~$n_2$.
  
Finally, we are ready to derive that $n_2=0$ by estimating the total number of embeddings. We group them  depending on how many points of $V_2$ are used. As we observed before, there are $n_0 n_1^{\fallingfact{3}}+n_0^{\fallingfact{3}}n_1$ embeddings that do not use $V_2$ at all. Next consider the case that exactly one vertex $i$ of $H$ is mapped into $V_2$. Then $i\in \{1,3\}$ (because of the red edge $02$). If, say $i=1$ then $0,3$ have to go to the same part $V'_i$ and $2$ has to go to the other part $V'_{1-i}$. Thus, by the symmetry between $1,3\in V(H)$, there are at most $2n_2\left(n_0 n_1^{\fallingfact{2}}+n_0^{\fallingfact{2}}n_1\right)$ embeddings that use exactly one vertex of~$V_2$. Also, there are at most $\binom{n_2}2\cdot (n-n_2)^2+o(n^2)$ embeddings with exactly two points in $V_2$. Indeed,  once we chose an unordered pair $uw\in\binom{V_2}2$, its adjacency uniquely determines which pair in $V(H)$ is mapped to it (namely, it is $\{0,2\}$ if $uw\in E(G)$ and $\{1,3\}$ otherwise) and then there are at most $(1/2+o(1)) \binom{n-n_2}2$ choices of a suitable unordered pair outside of $V_2$, and then there are $4$ ways to order these two unordered pairs. Finally, the number of embeddings with at least three vertices from $V_2$ is at most $O(n_2^3n)=O(n)$.

Thus, if we denote $c:=(n/2-n_0)/\sqrt{n}$ (and substitute $n_0=n/2+c\sqrt{n}$ and $n_1=n/2-c\sqrt{n}-n_2$) and use that $c=o(\sqrt{n})$ and $n_2=O(1)$, we get that
     \begin{eqnarray}
     0&\le & n_0 n_1^{\fallingfact{3}}+n_0^{\fallingfact{3}}n_1+ 2n_2\left(n_0 n_1^{\fallingfact{2}}+n_0^{\fallingfact{2}}n_1\right)+\binom{n_2}2\cdot (n-n_2)^2-\frac18\,n^{\fallingfact{4}}-\frac{3}{4}\, n^2+o(n^2)\nonumber\\
     &=& 
     p(c,n_2)
     n^2
        +\left(-4c^3n_2-2cn_2^2\right)n^{3/2}
        +o(n^2),\label{eq:H6Temp}
     \end{eqnarray}
     where we denote $p(c,n_2):=-2 c^4+3c^2- 2 c^2 n_2-\frac{n_2^2}{4}+\frac{3n_2}{4}-\frac{9}{8}$.
It follows that $c=O(1)$; thus the second term in~\eqref{eq:H6Temp} (that is, the one at $n^{3/2}$) can be absorbed into the error term. We conclude that $0\le p(c,n_2)+o(1)$. 
If we fix $n_2$ and vary $c$, then the derivative of the (even degree-4) polynomial of $p(c,n_2)$ has three roots: $c=0$ and $\frac12 \pm\sqrt{3-2n_2}$. Thus, if $n_2\ge 2$ then the only real root is  $c=0$ and $p(c,n_2)\le p(0,n_2)=-{n_2^2}/{4}+{3n_2}/{4}-{9}/{8}$, which is at most $p(0,3/2)=-9/16$, a contradiction. Finally, if $n_2=1$ then $p(c,1)=-2c^4+c^2-5/8$ is at most $p(1/2,1)=-1/2$, a contradiction again. Thus $n_2=0$, proving Theorem~\ref{th:H6Exact}.\epf

\subsection{Graph $H_5$}\label{se:H5}

Recall that $H_5$ is the red $3$-edge star $\MultiPartite{3,\,1}$ with one blue edge inserted. As we mentioned in the introduction, we have not been able to round the floating-point matrices returned by the plain flag algebra method for the $\obj[H_{5}]{}$-problem. However, we will show here (in Lemma~\ref{lm:CompletePartite}) that in order to solve the semi-inducibility problem for $H_5$  it is enough to work with complete partite graphs only. 
Specifically, in Section~\ref{se:SymmTheory} we present a sufficient criterion from \LPSS{} for a version of perfect stability to hold.  
Then, in Section~\ref{se:SuffSymm}, we show that it applies to the semi-inducibility problem for~$H_5$. The computational part (when we solve $H_5$-semi-inducibility problem in the theory of complete partite graphs) is discussed in Section~\ref{se:H5CP}.

\subsubsection{Symmetrizable problems}\label{se:SymmTheory}

First, we need to give some further definitions.  Given $\kappa\ge 2$ and $\gamma:\C F_\kappa\to\I R$, let the parameters $\obj[\gamma]{G}$, $\obj[\gamma]{n}$ and $ \obj[\gamma]{}$ be as defined in Section~\ref{se:gamma}.

For a graph $G$ and distinct vertices $u,w\in V(G)$, let the graph $G_{uw}$ be obtained from $G$ by making $w$ a clone of~$u$ (thus the new neighbourhood of $w$ is $\Gamma_G(u)\setminus\{w\}$). 
We call the $\obj[\gamma]{}$-problem \emph{strongly symmetrisable} if, for every graph $G$ and every non-edge $uw\in E(\O G)$, it holds that
 \beq{eq:symm}
 2\obj[\gamma]{G}\le \obj[\gamma]{G_{wu}}+\obj[\gamma]{G_{uw}}.
 \eeq

It is easy to show that, if the problem is strongly symmetrisable, then it is symmetrisable, as defined in \LPSS{Definition 1}. Since the latter definition is rather long (in order to apply for a wider class of problem), we do not state it here. Instead, we verify, for example, the first  part of \LPSS{Definition 1}, namely that one can transform an arbitrary graph $G_0$ on $[n]$ into a complete partite graph in steps so that each step changes only $o(n^2)$ pairs (in fact, it will be at most $n-1$ pairs here) and does not decrease the objective function. Inductively on $i$, we construct a graph $G_i$ whose restriction to $[i-1]$ is complete partite. Suppose that $i\ge 1$ and already have $G_{i-1}$. If $i$ is connected to every vertex of $[i]$ in $G_{i-1}$ then we take $G_{i}:=G_{i-1}$ (with $\{i\}$ forming a new separate part), so suppose otherwise. Take any $j\in [i]$ non-adjacent to $i$ in $G_{i-1}$. 
If $\obj[\gamma]{(G_{i-1})_{ji}}\ge \obj[\gamma]{G_{i-1}}$ then we let $G_i:=(G_{i-1})_{ji}$ (where we replace $i$ by a clone of $j$ in $G_{i-1}$). So suppose that $\obj[\gamma]{(G_{i-1})_{ji}}< \obj[\gamma]{G_{i-1}}$. By~\eqref{eq:symm} we have that $\obj[\gamma]{(G_{i-1})_{ij}}> \obj[\gamma]{G_{i-1}}$. Let $J=\{j_0,\dots,j_{m-1}\}$ denote the set of $G_{i-1}$-clones of $j$ that lie inside $[i]$. 
We get $G_i$ in $m$ steps where in the $t$-th step we symmetrise $j_t$ into $i$ in the current graph. By~\eqref{eq:symm} and induction on $t$, the $t$-th step increases the objective function (because the opposite symmetrisation decreases it), finishing the proof. A direct consequence is that, for every integer $n\ge \kappa$, at least one extremal graph attaining the value of $\obj[\gamma]{n}$ is complete partite.

Let $\C P$ denote the set of all complete partite graphs (up to isomorphism) and let $\O{\C P}$ be the set of sequences $(x_i)_{i=1}^\infty$ such that $x_1\ge x_2\ge \ldots$, each $x_i$ is non-negative and $\sum_{i=1}^\infty x_i\le 1$. 
There is a natural map that sends a complete $m$-partite graph $\MultiPartite{n_1,\,\dots,\,n_{m}}$ with $n_1\ge\dots\ge n_m$ to the sequence $(n_1/n,\ldots,n_m/n,0,\ldots)\in\O{\C P}$, where $n:=n_1+\dots+n_m$. 
The set  $\O{\C P}$ with the pointwise convergence can be viewed as the limit space for complete partite graphs and its elements can be identified with graphons of a very special form (namely, those for which the induced density of the 3-vertex 1-edge graph $\O P_3$ is zero). 
We refer the reader to \LPSS{} for further details. Observe that $x_0:=1-\sum_{i=1}^\infty x_i$ may be positive in the limit; combinatorially, it refers for the total proportion of vertices in parts whose relative size tends to 0. Also, the following definition gives, for every element $\V x\in \O{\C P}$ a sequence of complete partite graphs that converges to~$\V x$.

 \begin{definition}
 \label{def:Gnx}
 Given $n\in\mathbb{N}$ and $\V x = (x_1,x_2,\ldots)$ with $x_1 \geq x_2 \geq \ldots \geq 0$ and $x_0 := 1-\sum_{i \geq 1}x_i \geq 0$,
 define a complete partite graph $G_{n,\V x}$ with vertex set $[n]$, parts $V_1,\ldots,V_m$ for some $m$ and a set $V_0$ of universal vertices, i.e., $|V_0|$ singleton parts, as follows.
 If $x_0=0$, take a partition $[n]=V_1 \cup \ldots \cup V_m$ with $\left|\, |V_i|-x_i n\,\right| < 1$
 and let $V_0=\emptyset$.
 Otherwise, for all $i \geq 1$ with $x_i n \geq 2$, let $|V_i|=\lfloor x_i n \rfloor$ and let $V_0$ consist of the remaining vertices in $[n]$.
 \end{definition}

Next, we define a version of perfect stability from~\LPSS{}.

 \newcommand{\CP}{\mbox{$\C P$}}
 \begin{definition}[Perfect \CP-Stability]
 \label{de:CPPerfectStab}
 The $\obj[\gamma]{}$-problem is \emph{perfectly \CP-stable} if there exists $C>0$ such that for every graph $G$ of order $n\geq C$ there is a complete partite graph $G'$ of the same order $n$
 such that
 $$\dedit(G,G')\leq C(\obj[\gamma]{n}-\obj[\gamma]{G})\,n^2.$$
 \end{definition}
 
Note that this property, although very closely related to our definition of perfect stability, is slightly different since the number of parts of $G'$ is not specified. 

Before we can state the sufficient condition for \CP-stability from \LPSS{}, we need to give some further definitions.  
Define
 \begin{equation}\label{eq:CPlambda}
     \obj[\gamma]{\V x} := \lim_{n\to \infty}\obj[\gamma]{G_{n,\V x}}
 \end{equation}
and let $\OPT(\obj[\gamma]{})$ be the set of $\V x\in\C P$ that maximise this function.
 One can show that the limit in~\eqref{eq:CPlambda} exists and write an explicit formula for $\obj[\gamma]{\V x}$ (which is a polynomial of $\V x$ if $\V x$ has finite support and an absolutely convergent power series otherwise). 
 Also, by the compactness of $\O{\C P}$ and the continuity of the function $\obj[\gamma]{}:\O{\C P}\to\I R$, the set $\OPT(\obj[\gamma]{})$ is non-empty.

In the notation of Definition~\ref{def:Gnx}, if $G'$ is a graph obtained by adding a new vertex $z$ to $G=G_{n, \V x}$, we say that $z$ is a \emph{\CP-clone} of $u \in V(G)$ if $u \in V_0$ and $\Gamma_{G'}(z)=V(G)$, or if $u \notin V_0$ and $\Gamma_{G'}(z)=\Gamma_G(u)$.

 \begin{theorem}[\LPSS{Theorem 1.1}]\label{thm-main-rough}
 Let $\obj[\gamma]{}$ be a symmetrisable function such that the set $\OPT(\obj[\gamma]{})$ has finitely many elements.
 Suppose also that there exists $c>0$ such that the following hold for every large $n$ and every maximiser $\V x \in \OPT(\obj[\gamma]{})$, where $G:=G_{n,\V x}$.
 \begin{enumerate}[(i)]
 \item\label{it:main-rough1} For all distinct $u,w \in V(G)$ we have $\obj[\gamma]{G}-\obj[\gamma]{G \oplus uw} \geq cn^{-2}$..
 \item\label{it:main-rough2}
 If $G_u$ is obtained from $G$ by adding a new vertex $u$ which is complete or empty to each part of $G$
 (where each $V_i$, $i \in [m]$ is a part and we have $|V_0|$ singleton parts)
 then the minimum number of edits at $u$ needed to make $u$ a \CP-clone of some existing vertex of $G$ is at most $n(\obj[\gamma]{G}-\obj[\gamma]{G_u, u})/c$.
 \end{enumerate}
 Then the $\obj[\gamma]{}$-problem is perfectly \CP-stable.
 \end{theorem}

Note that if every $\V x\in \OPT(\obj[\gamma]{})$ has $m<\infty$ non-zero terms and satisfies $\sum_{i=1}^m x_i=1$ then Property~\ref{it:main-rough1} of Theorem~\ref{thm-main-rough} is the same as our $\obj[\gamma]{}$-flip-averseness of $K_m$ while Property~\ref{it:main-rough2} can be shown to be implied by the $\obj[\gamma]{}$-strictness of $K_m$. (Note that when we attach a new vertex $u$ in Property~\ref{it:main-rough2}, we consider only those attachments where, for each part $V_i$, the pairs between $u$ and $V_i$ are all edges or all non-edges;
Property~\ref{it:main-rough2} relates the contribution of $u$ to the number of edits so that it handles correctly the general case in \LPSS{} when $\V x$ may have infinitely many positive entries.)

\subsubsection{Sufficient condition for a semi-inducibility problem to be symmetrisable}\label{se:SuffSymm}

The following lemma gives a sufficient condition for the $H$-semi-inducibility problem to be symmetrisable. It generalises the analogous inducibility result of Brown and Sidorenko~\cite{BrownSidorenko94} which is a special case when $H$ is a 2-colouring of a complete graph. It is proved by an easy adaptation of the argument from~\cite{BrownSidorenko94}.

\begin{lemma}\label{lm:CompletePartite} Let $H$ be a blue-red edge-coloured graph 
such that the red graph is the complete partite with parts $R_0\cup\dots\cup R_{m-1}$ and, for  every $i\in [m]$ and every two vertices $u,w\in R_i$, the blue neighbourhoods of $u$ and $w$ in $R_i\setminus\{u,w\}$ are nested (that is, one is a subset of the other, allowing equality). Then the $H$-semi-inducibility problem is strongly symmetrisable.\end{lemma}

\bpf For a graph $G$, let $\mathrm{Hom}(H,G)$ be the set of embeddings of $H$ into~$G$. Recall that these are injections $V(H)\to V(G)$ that map red (resp.\ blue) edges to edges (resp.\ non-edges) of $G$. Thus $|\mathrm{Hom}(H,G)|=\Obj[H]{G}$. 

Take any graph $G$ and $xy\in E(\O G)$. In order to prove~\eqref{eq:symm} it is enough to prove that, for every distinct $u,w\in V(H)$ and any injection $g:V(H)\setminus\{u,w\}\to V(G)\setminus \{x,y\}$, we have
 \beq{eq:symm2}
 2\,\left|\mathrm{Hom}(H,G)\cap \mathrm{Ext}(g)\right| \le \left|\mathrm{Hom}(H,G_{xy})\cap \mathrm{Ext}(g)\right|+\left|\mathrm{Hom}(H,G_{yx})\cap \mathrm{Ext}(g)\right|,
 \eeq
 where $\mathrm{Ext}(g)$ consists of all \emph{extensions} $f$ of $g$, that is, injective maps $V(H)\to V(G)$ such that $f(v)=g(v)$ for every $v\in V(H)\setminus\{u,w\}$.
 Indeed, if we sum~\eqref{eq:symm2} over all unordered pairs $uw\in\binom{V(H)}2$ and all injections $g:V(H)\setminus\{u,w\}\to V(G)\setminus \{x,y\}$ then every embedding $f$ from $H$ into any one of $G$, $G_{xy}$ or $G_{yx}$ is counted exactly $\binom{v(H)}2$ times.
 
Take any extension $f$ of $g$. Let $\phi$ be the identity map on $V(G)$, except it swaps $x$ and $y$. Note that $\phi$  is an involution, that is, $\phi\circ \phi$ is the identity map.

If the image of $f$ is disjoint from $\{x,y\}$ and  $f$ is an embedding of $H$ into $G$ then it is trivially an embedding of $H$  into each of $G_{xy}$ and $G_{yx}$, so it contributes $2$ to both sides of~\eqref{eq:symm2}, as desired.  

Next, let us consider those extensions $f\in\mathrm{Ext}(g)$ such that the image of $f$ contains exactly one of $x$ and $y$. If $f$ is an embedding of $H$ into $G$ and, say, $f(u)=x$ (and thus $f(w)\not=y$) then both $f$ and the composition $\phi\circ f$ are two different embeddings into $G_{xy}$ (since $y$ is a clone of $x$ in $G_{xy}$), contributing $2$ to the right-hand side of~\eqref{eq:symm2}. We do not overcount in the right-hand side: if $\phi\circ f$ also happens to be an embedding into $G$, then the contribution of $\phi\circ f$ is also counted for $G_{yx}$ (and the pair $\{f,\phi\circ f\}$ contributes 4 to both sides of~\eqref{eq:symm2}).

Finally, it remain to look at extensions $f$ whose image contains both $x$ and~$y$. Since $\phi$ is an involution, all such maps are partitioned into pairs $\{f,\phi\circ f\}$. Consider the contribution of any such pair.
There is nothing to prove if none of $f$ and $\phi\circ f$ is an embedding of $H$ into $G$.
Suppose that both $f$ and $\phi\circ f$ are embeddings of $H$ into $G$. By $xy\in E(\O G)$, the pair $uw$ is not a red edge of~$H$. Thus $u$ and $w$ are in the same part $R_i$ and have the same red neighbourhoods in~$H$. Also, for  every vertex $v\in V(H)\setminus\{u,w\}$, which is a blue neighbour of at least one of $u$ and $w$, we have that $f(v)\in \Gamma_{\O G}(x)\cap \Gamma_{\O G}(y)$. Thus both $f$ and $\phi\circ f$ are embeddings into both $G_{xy}$ and $G_{yx}$, and  pair $\{f,\phi\circ f\}$ contributes $4$ to both sides of~\eqref{eq:symm2}. Finally, suppose that $f$ is an embedding of $H$ into $G$ but $\phi\circ f$ is not. By the former property and by $xy\in E(\O G)$, there is $i\in [m]$ such that $u,w\in R_i$. Thus $u$ and $w$ have the same red neighbourhoods in $H$ and nested blue neighbourhoods in $R_i\setminus\{u,w\}$, say that of $u$ being a subset of that of~$w$. Thus the constraints on the possible image of $f(u)$ (for $f\in \mathrm{Ext}(g)$) are  weaker than those for $f(w)$. We conclude that both $f$ and $\phi\circ f$ are embeddings into $G_{f(w)f(u)}$. Thus the contribution of $\{f,\phi\circ f\}$ to the right-hand size of~\eqref{eq:symm2} is at least~2, as desired. This finishes the proof of the lemma.\epf

\subsubsection{Calculations for $H_5$}\label{se:H5CP}

Clearly, Lemma~\ref{lm:CompletePartite} applies when $H=H_5$
is the red 3-star with one blue edge inserted. Thus, in order to determine the semi-inducibility constant of $H$, it is enough to solve the problem for complete partite graphs.
While one should in principle be able to determine the optimal number of parts (5 in this case) and the optimal limiting ratios (the vector $(\alpha,\dots,\alpha,(1-4\alpha))$ with $\alpha=(13-\sqrt{57})/{56}$) by hand, this task seems still quite messy. Instead, we recourse to flag algebras.

\begin{theorem}\label{th:semiind5CP}
Let \[H := H_5 = \DrawSemiGraphLabeled{4}{0/1,0/2,0/3}{2/3}{0,2,1,3}\] be the graph on $[4]$ with the red edge set $\{01, 02 ,03\}$ and the blue set $\{13\}$. Then the maximum value of $\obj[H]{G}$ over complete partite graphs $G$ of order $n\to\infty$ is $\frac{-171 \sqrt{57} + 7879}{43904}+o(1)$. Moreover, there is $C$ such that for every complete partite graph $G$ of order $n\ge 1/C$ there is a complete $5$-partite graph $G'$ with 
\begin{equation}\label{eq:H5CPStabilty}
\dedit(G,G')\leq C(\obj[\gamma]{n}-\obj[\gamma]{G})\,n^2.
\end{equation}
\end{theorem}

\bpf 
The first part of the theorem is proved by the standard flag algebra approach with $\N=8$ within the theory of complete partite graphs. This is a hereditary family of graphs (characterised by forbidding $\O{P}_3$ as an induced subgraph) so flag algebras can be applied to it.

Let us show how we derive the second part of the theorem from the obtained certificate.

Our certificate 
satisfies that the slack $c_F$ is positive for
every 8-vertex complete partite graph $F$ with more than 5 parts. Hence up to removing $o(n)$ vertices, every almost extremal complete partite graph is 5-partite.
Also, it satisfies \PST{Theorem 7.1}, including its Condition~(i), with $\tau=K_4$. (The statement of \PST{Theorem 7.1} allows for any hereditary family of graphs which is closed under taking blowups, in particular it applies to~$\C P$.)
Thus \PST{Theorem 7.1} implies in particular that there is the unique (up to an automorphism of $K_5$, that is, any permutation of indices) vector
$\V x\in\I S_5$ that maximises $\obj[H]{K_5(\V x)}$. Thus $\V a:=(\alpha,\alpha,\alpha,\alpha,1-4\alpha)$ is the unique, up to a permutation, $(\obj[H_5]{},K_5)$-optimal vector.

The claim about the existence of $C$ is another conclusion of \PST{Theorem 7.1}, namely that the perfect $K_5$-stability property (defined in the obvious way) holds within  the theory of complete partite graphs.\epf

Combining Lemma~\ref{lm:CompletePartite} and Theorem~\ref{th:semiind5CP}, we can obtain the upper bound (and thus the value of $\obj[H_5]{}$) stated in Table~\ref{ta:semiind}: 

\begin{theorem}\label{th:semiind5}
 It holds that $\obj[H_5]{}=\frac{-171 \sqrt{57} + 7879}{43904}$ and the semi-inducibility problem for $H_5$ is perfectly $K_5$-stable.
 \end{theorem}
 \bpf This extremal problem is symmetrisable by  Lemma~\ref{lm:CompletePartite}. Furthermore, Theorem~\ref{th:semiind5CP} gives that the complete partite version of the problem has unique maximiser in the limit space $\O{\C P}$, namely the sequence made from $\V a:=(\alpha,\alpha,\alpha,\alpha,1-4\alpha)$ by attaching zeros. This corresponds to taking blowups of $B:=K_5$.

Next, we would like to check  the flip averseness and Condition~\ref{it:main-rough2} of Theorem~\ref{thm-main-rough} for $\obj[H_5]{}$. While this task is straightforward, the number of cases is considerable: up to symmetries we have to check 4 possible flip pairs (since there are 2 different part ratios) and 10 possible attachments  (namely, the new vertex can be fully attached to any number $i\in [5]$ of parts of ratio $\alpha$ times the two options of being full or empty to the remaining part of ratio $1-4\alpha$). So we delegate this task to computer and provide scripts that verify both flip averseness and the strictness using exact arithmetic. Recall that strictness implies  Condition~\ref{it:main-rough2} of Theorem~\ref{thm-main-rough}. 
 Thus Theorem~\ref{thm-main-rough} applies here and gives that the semi-inducibility problem for $H_5$ is perfectly \CP-stable. This, when combined with~\eqref{eq:H5CPStabilty}, gives that the problem is in fact perfectly $K_5$-stable.\epf

\section{Concluding remarks}

Let us remark that, although we proved that the semi-inducibility problem for $H_4$ is perfectly $K_3$-stable in Section~\ref{se:H4}, our flag algebra certificate with $\N=4$ does not satisfy the sufficient condition of \PST{Theorem~7.1}. (Namely  the 4-vertex clique $K_4$ does not admit a homomorphism into the base graph $K_3$ but has zero slack.) However, the information contained in the certificate is enough to give alternative proofs of Claims~\ref{cl:H4OP3} and~\ref{cl:3Parts}. For example, every $4$-vertex graph $F$ which is not complete partite has positive slack $c_F$ in~\eqref{eq:FAMain}, which implies that the density of $F$ in every almost extremal graph $G$ is~$o(1)$. Since (every 4-vertex extension of) $\O{P}_3$ is not complete partite,  Claim~\ref{cl:H4OP3} (which is equivalent to the induced density of $\O{P}_3$ being $o(1)$) follows.

Flag algebras worked remarkably well for the semi-inducibility of non-complete 4-vertex graphs. Recall that if we are interested only in the value of the semi-inducibility constant $\obj[H_i]{}$, $0\le i<18$, then the sharp upper bounds can be determined by this method for every $i$ except perhaps $i=3$ and $i=5$. (In fact, it is possible that the method also works for these two cases but we were just unable to round the floating-point matrices returned by the solver.) One informal explanation for such efficiency is that, since $\{2,3\}$ is not an edge of $H_i$, we can compute $\Obj[H_i]{G}$ up to error $O(n^3)$ by summing over all suitable $f(0),f(1)\in V(G)$ of the product of the number of choices for $f(2)$ and $f(3)$. This can be viewed as some scalar product of the  vectors of 3-vertex subgraph densities in $G$ rooted on $(f(0),f(1))$. Since positive semi-definite matrices appear naturally for problems involving scalar products (e.g.\ as Gram matrices), this gives some intuition why the flag algebra method is so successful here. 

Using flag algebras, the first author has obtained a number of exact results on the semi-inducibility constant for 5-vertex graphs. These findings will appear in a subsequent paper.

We should also mention that the computer-generated certificates  obtained early in the project were quite helpful in finding matching lower bounds as well as in proving in a computer-free way some of the results (once the authors were confident of their validity). This approach may be useful for other problems. The first author made publicly available not only the verifier code but also the code that was used for generating the certificates. It should  not be hard to adapt it to apply flag algebra calculations to other suitable  extremal problems. We hope that this package will be generally useful to the research community.

\section*{Acknowledgements}

We would like to thank the anonymous reviewers for their constructive feedback and insightful suggestions. Both authors were supported by ERC Advanced Grant 101020255.

\section*{Declaration on the use of AI}

Claude Opus 4.7 Adaptive Thinking was used to assist in proofreading. All mathematical arguments, results, and conclusions were written and verified by the authors, or by a deterministic computer program included in the ancillary files.

\bibliography{bibexport}

\end{document}